\documentclass[12pt,a4paper]{article}
\usepackage{latexsym,amssymb,amsthm}
\RequirePackage[dvips,dvipdfm]{graphicx}

\title{\bf The structure of decomposition of a triconnected graph}

\author{D.\,V.\,Karpov \and A.\,V.\,Pastor}

\date{}

\begin{document}

\def\q#1.{{\bf #1.}}
\def\P{{\rm Part}}
\def\I{{\rm Int}}
\def\R{{\rm Bound}}
\def\O{{\rm Nb}}

\def\D{{\rm Dep}}
\def\Cl{{\rm Class}}
\def\Co{{\rm Comp}}
\def\Str{{\rm Struct}}

\sloppy
\righthyphenmin=2
\exhyphenpenalty=10000
\binoppenalty=8000
\relpenalty=8000

\renewcommand*{\proofname}{\bf Proof}
\newtheorem{thm}{Theorem}
\newtheorem{lem}{Lemma}
\newtheorem{cor}{Corollary}
\theoremstyle{definition}
\newtheorem{defin}{Definition}
\theoremstyle{remark}
\newtheorem{rem}{\bf Remark}
\newcommand{\ov}{\overline}

\maketitle

\section*{Introduction}

The structure of decomposition of a connected graph by its cutpoints 
(i.e., vertices which deleting makes graph disconnected) is well known~\cite{O,X}. It is convenient
to describe this structure with the help of so-called  {\it tree of blocks and cutpoints}. The vertices of this tree are  cutpoints and  parts of decomposition of the graph by its cutpoints.

In  1966  W.\,T.\,Tutte~\cite{T} described the structure of relative disposition of 2-vertex cutsets in biconnected graphs and  showed, that this structure has much in common with the the structure of cutpoints. 
In particular, a construction of the tree of blocks for biconnected graph was introduced. Analogous constructions 
were considered also in the works~\cite{Mac,HT}.

Attempts of development of analogous constructions for graphs of greater connectivity were done in the works~\cite{Hoh,KP,k02}. But  significant difficulties appear in this process. They are concerned with 
the fact that two $k$-vertex cutsets of a $k$-connected graph can be dependent, i.e. after deleting of one of them, vertices  of another could appear in different connected components.
This leads to non-uniqueness of resulting constructions of a tree of blocks for $k$-connected graphs.
These constructions are essentially  dependent on the order in which cutsets were chosen during the process. Moreover, such constructions take account of not all  $k$-vertex cutsets:
splitting the graph by one cutset we automatically loose information about all cutsets dependent with the chosen one. In the works~\cite{KP,k02} these difficulties were partly overcome for graphs, satisfying
some additional condition. But in general case the question of how to describe the structure of decomposition
of a $k$-connected graph by its $k$-vertex cutsets for $k\ge3$ remained open.

In the work~\cite{k05} it was developed a new method of studying of the structure  of relative disposition of  $k$-vertex cutsets of a  $k$-connected graph~--- the theorem of decomposition. With the help of this theorem
several results for the case of arbitrary $k$ were obtained. As an illustration of work of new method in the end of the work~\cite{k05} one can see rather simple and visual description of the structure of 
2-vertex cutsets of a biconnected graph. In general, this description is similar to the construction of Tutte~\cite{T}, but it is a good illustration of efficiency of the new method.

This paper is devoted to studying of the structure of relative disposition of 3-vertex cutsets in a (vertex)
triconnected graph. We will use the theorem of decomposition and as a  result we obtain
a description, similar to  analogous structure of a biconnected graph~\cite{T,k05}.

\section{Basic notations}

Always in our paper we consider simple undirected  finite graphs without loops and multiple edges.

We use the following notations and definitions. 
For a graph $G$ we denote the set of its vertices by  $V(G)$ and the set of its edges by $E(G)$.
We denote the {\it degree} of a vertex $x$ in the graph $G$ by $d_G(x)$.

We call two vertices {\it connected}, if there is a path between them.
By a {\it connected component} of a graph we always mean its maximal  (with respect to inclusion) set of pairwise 
connected vertices.

A graph~$G$ is called  {\it $k$-connected}, if it contains at least $k+1$ vertices and remains connected after deleting  any   $k-1$ vertices.
In particular, for  $k=2$ such a graph is called {\it biconnected}, and for $k=3$~--- {\it triconnected}.

For any set of edges  $E\subset E(G)$ we, as usual, denote by $G-E$ the graph obtained from~$G$  after deletion 
of edges of the set~$E$.  For $e\in E(G)$ we set $G-e=G-\{e\}$.

For any set of vertices $V\subset V(G)$ we denote by  $G-V$ the graph obtained from~$G$  after deletion 
of vertices of the set ~$V$ and all edges incident to deleted vertices.  
For  $v\in V(G)$ we set $G-v$=$G-\{v\}$.

For any set $M\subset V(G)\cup E(G)$ we denote by  $G-M$ the graph obtained from~$G$  after deletion
of all vertices and edges of the set~$M$ and all edges incident to deleted vertices.

During all our work let $G$ be a triconnected graph with $|V(G)|>6$. 

A set $S\subset V(G)$ is called a {\it cutset}, if the graph $G-S$ is disconnected.   
We denote the set of all cutsets of the graph~$G$  by~$\mathfrak R(G)$,
and the set of all  $3$-vertex cutsets of $G$ (we will call them simply 3-{\it cutsets})
--- by~$\mathfrak R_3(G)$.

We use   terminology of the work~\cite{k05}. We rewrite   definitions from~\cite{k05}, that we need, in the form convenient for triconnected graphs.

\begin{defin}
$1)$ Let  $R,X\subset V(G)$.
We say that~$R$ {\it splits}~$X$, if not all  vertices of the set $X\setminus R$ 
are  in the same connected component  of the graph ${G-R}$.

$2)$ Let $U,W \subset V(G)$.
We say that~$R$ {\it separates} a set~$U$ from a set~$W$, if
 $U\not\subset R$, $W\not\subset R$ and any two vertices $u\in U\setminus R$ and $w\in W\setminus R$ 
lie in different connected components of the graph~${G-R}$.

In the case  $U=\{u\}$ we say that~$R$ separate a vertex~$u$ from a set~$W$.
If $U=\{u\}$ and $W=\{w\}$, we say that~$R$ separates a vertex~$u$  from a vertex~$w$.
\end{defin}

\begin{defin}
1) We call sets $S,T\in \mathfrak R_3(G)$ {\it independent}, if~$S$ 
does not split~$T$ and~$T$ does not split~$S$. Otherwise, we call these sets 
{\it dependent}.

2) We assign to each set  $\mathfrak S\subset \mathfrak R_3(G)$ the 
{\it dependence graph} $\D(\mathfrak S)$, which vertices are  cutsets of~$\mathfrak S$, and two vertices are adjacent if and only if  correspondent cutsets are dependent.
\end{defin}

Thus, any set~$\mathfrak S$ is divided into  {\it dependence components}~--- 
subsets, correspondent to connected components of the graph~$\D(\mathfrak S)$.

It is easy to prove, that if~$T$ does not split~$S$, then~$S$ does not split~$T$, i.e. these cutsets are
independent (see~\cite{Hoh,KP}).

\begin{defin}
Let $\mathfrak S\subset {\mathfrak R_3}(G)$.

1)  A {\it part} of decomposition of the graph~$G$ by the set~$\mathfrak S$ 
(or a part of $\mathfrak S$-decomposition)~is a maximal (with respect
to inclusion) set $A\subset V(G)$ such, that  no cutset
$S\in \mathfrak S$ splits~$A$.
We denote by $\P(\mathfrak S)$ the set of all such parts.
If~$\mathfrak S$ consists of one cutset~$S$, then we denote
the set of all parts of  $\{S\}$-decomposition by $\P(S)$.

2)  Let $A\in \P(\mathfrak S)$. We call  {\it inner vertices} all vertices of~$A$ which do not belong to any cutset of~$\mathfrak S$. The set of all inner vertices of the part $A$  we call
  {\it interior} of the part~$A$ and denote by~$\I(A)$.  

We call {\it boundary vertices} all vertices of the part~$A$ belonging to any cutsets from~$\mathfrak S$.
The set of all such vertices  we call {\it boundary} of the part~$A$ and denote by $\R(A)$.

3) We call a part~$A$ {\it empty}, if $\I(A)=\varnothing$ and {\it nonempty} otherwise.
We call a part~$A$ {\it small}, if $|A|<3$ and {\it normal}, if $|A|\ge 3$.
\end{defin}

It is easy to see, that two different parts~$A_1,A_2\in \P(\mathfrak S)$ either have no common vertices, or $A_1\cap A_2$ is a subset of one of cutsets of~$\mathfrak S$.
It is proved in~\cite[theorem~2]{k05}, that~$\R(A)$  consists of all vertices of a part~$A$, which are adjacent to vertices outside~$A$ and~$\R(A)$ separates $\I(A)$ from $V(G)\setminus A$.

An important particular case of decomposition of triconnected graph by a set of 3-cutsets is
a decomposition by one 3-cutset~$S$.
It is clear, that  for any part  $F\in \P(S)$ its interior $\I(F)$ is a connected component 
of the graph $G-S$.

Since no subset of the cutset~$S$ is a cutset of the graph~$G$, then any vertex of~$S$ is adjacent 
to at least one vertex of~$\I(F)$,  hence, the induced subgraph of the graph~$G$ on the vertex set~$F$ is connected.

Note, that any vertex~$x$ of the graph~$G$ is adjacent to at least one other vertex~$y$.
Obviously, no cutset can separate~$x$ from~$y$, thus, for any set~$\mathfrak S\subset \mathfrak R_3(G)$ any  part~${A\in\P(\mathfrak S)}$ contains at least two vertices. 
Hence, any small part contains exactly two vertices.

\subsection{Dependent and independent cutsets}

Let $S,T\in\mathfrak R_3(G)$. 
Clearly these cutsets are independent if and only if there exists a part $F\in\P(S)$ which contains~$T$. It was proved in~\cite[lemma~1]{k05} that if there exists a part~$A\in \P(S)$ such that $\I(A)\cap T=\varnothing$, then the cutsets~$S$ and~$T$ are independent. 
Decomposition of the graph by a pair of dependent cutsets is described in the following lemma.

\begin{lem}[{\cite[lemma~7]{k05}}]
\label{lds0}
Let $G$ be a $k$-connected graph and cutsets~$S,T\in {\mathfrak R}_k(G)$ be dependent. Let $\P(S)=\{F_1,\ldots,F_n\}$ and 
$\P(T)=\{H_1,\ldots,H_m\}$.
For all  ${i\in\{1,\ldots,n\}}$ and $j\in\{1,\ldots,m\}$ we set 
   $$P=S\cap T, \quad S_j=S\cap \I(H_j), \quad T_i=T\cap \I(F_i), \quad G_{i,j}=F_i\cap H_j.$$
Then 
   $$\P(\{S,T\})=\{G_{i,j}\}_{i\in\{1,\ldots,n\},\:j\in\{1,\ldots,m\}}, \qquad 
     \R(G_{i,j})=P\cup T_i\cup S_j,$$
moreover, $T_i\ne\varnothing$ for all $i\in\{1,\ldots,n\}$ and~$S_j\ne\varnothing$ for all  $j\in\{1,\ldots,m\}$.
\end{lem}

The statement of lemma~\ref{lds0} is correct for $k$-cutsets of a  $k$-connected graph for all~$k$. In the case $k=3$, which is interesting to us, it is easy to derive the following statements from this lemma. 
(Notations are the same as in the lemma).

\begin{cor}
\label{l1c1}
Let cutsets $S,T\in\mathfrak R_3(G)$ be dependent. Then $|{S\cap T}|\le1$, and
each of these cutsets splits $G$ into not more than~$3$ parts.

\begin{proof}
It is easy to see, that $m,n\le3$, since all the sets~$T_i$ and~$S_j$ are nonempty.
Obviously, $m,n\ge2$, hence, $|P|\le1$.
\end{proof}
\end{cor}

\begin{cor}
\label{l1c2} 
If $S\cap T=\varnothing$, then~$\P(\{S,T\})$ contains at least one small part. 
 $|G_{i,j}|=2$ if and only if $|T_i|=|S_j|=1$. 
Any small part  $G_{i,j}\in\P(\{S,T\})$ consists of two vertices~$u$ and~$v$, where  $u\in T$ 
and $v\in S$. Moreover,  $v$ is the only vertex of the part~$H_j$ adjacent to~$u$, and
$u$ is the only vertex of the part~$F_i$ adjacent to~$v$.

\begin{proof}
Since $|T_1|+|T_2|\le3$, at least one of the sets~$T_1$ or~$T_2$ consists of one vertex.   	
Without loss of generality we may assume, that it is~$T_1$. Analogously, we may assume, that~$|S_1|=1$. 
Then $|\R(G_{1,1})|={|T_1\cup S_1|=2}$. Hence, $\I(G_{1,1})=\varnothing$ and $|G_{1,1}|=2$. 
Similarly, each part~$G_{i,j}$ for which $|T_i|=|S_j|=1$ is small.

Let $|G_{i,j}|=2$. Obviously, then $|T_i|=|S_j|=1$ and $G_{i,j}=\{u,v\}$, where $T_i=\{u\}$, $S_j=\{v\}$. Since $u\in T$, then the vertex~$u$ must be adjacent to at least one inner vertex of the part~$H_j$. 
On the other side, since $u\in \I(F_i)$, then the vertex~$u$ can be adjacent only to  vertices
of the part~$F_i$. As well, $F_i\cap H_j=G_{i,j}=\{u,v\}$, hence, $v$ is the only vertex of the part~$H_j$ which can be adjacent to~$u$. Thus, the vertices~$u$ and~$v$ are adjacent and $v$ is the only vertex of the part~$H_j$ adjacent to~$u$. Analogously, $u$ is the only vertex of the part~$F_i$ adjacent to~$v$.
\end{proof}
\end{cor}

\begin{cor}
\label{l1c3}
If $|S\cap T|=1$, then  $m=n=2$ and  $\P(\{S,T\})$ contains no small parts. 
Any empty part of $\P(\{S,T\})$ consists of exactly three vertices~$u$, $v$ and~$p$,
where~$u\in S\setminus T$, $v\in T\setminus S$ and $P=\{p\}$. Moreover, the vertices~$u$  and~$v$ are adjacent.

\begin{proof}
Since all sets~$T_1$, $T_2$, $S_1$, $S_2$, $P$ are nonempty, we obtain  $m=n=2$ and
 $|T_1|=|T_2|=|S_1|=|S_2|=|P|=1$. Any part~$G_{i,j}$ must contain at least one  vertex from the sets~$T_i$, 
$S_j$ and~$P$, i.e., at least  3 vertices. Thus, there are no small parts.

If $\I(G_{i,j})=\varnothing$, then  $G_{i,j}=\R(G_{i,j})=T_i\cup S_j\cup P$.
Hence, $|G_{i,j}|=3$. Let $T_i=\{u\}$, $S_j=\{v\}$, $P=\{p\}$. Then
$G_{i,j}=\{u,v,p\}$. We can prove, that the vertices~$u$ and~$v$ are adjacent, as well as in corollary~\ref{l1c2}.
\end{proof}
\end{cor}

\begin{rem}
\label{r0}
We can   exclude the case $m=n=3$, because in this case $P=\varnothing$ and
$|T_1|=|T_2|=|T_3|=|S_1|=|S_2|=|S_3|=1$. Obviously, then all parts of $\P(\{S,T\})$ are small. Thus, $V(G)$
 consists of  vertices of the sets~$S$ and~$T$, i.e. $|V(G)|=6$.
As it was written in the beginning of our paper, we do not consider such graphs.
(It is easy to see, that in this case the graph~$G$ is isomorphic to~$K_{3,3}$.)
\end{rem}

\begin{lem}
\label{lds1}
Let sets $S,T\in\mathfrak R_3(G)$ and parts $A_1,\ldots,A_k\in \P(T)$ be such that
 $S\cap \I(A_i)=\varnothing$ for $i\in \{1,\ldots, k\}$. Then the set~$S$ does not 
split $A=\cup_{i=1}^k \:A_i$.

\begin{proof}
It is easy to see, that vertices of each of sets $\I(A_1),\ldots,\I(A_k)$ are connected in $G-S$. Every vertex of  a nonempty set~$T\setminus S$ is adjacent to a vertex of each set  $\I(A_1),\ldots,\I(A_k)$. Hence, 
$S$ does not  split~$A$.
\end{proof}
\end{lem}

\section{Basic structures}

In this section we describe basic structures, which dependent cutsets can form, and investigate basic properties
of these structures.

\subsection{Flowers in triconnected graphs}

Consider a tuple  $F=(p;q_1,\ldots,q_m)$ of vertices of our graph~$G$ (here $m\ge4$), in which the vertices $q_1,\ldots,q_m$ are  {\it cyclic ordered}. We will set, that a cyclic permutation of the set $q_1,\ldots,q_m$ does not change the tuple~$F$.
Let us introduce the notation $Q_{i,j}=\{q_i, q_j, p\}$. Let  $\mathfrak R(F)$ consist of sets
$Q_{i,j}$ for all pairs of different and non-neighboring in the cyclic order indexes~$i$ and~$j$.

\begin{defin}
We say, that a tuple $F=(p;q_1,\ldots,q_m)$ is a  {\it flower}, if there exists such a set  $\mathfrak S\subset \mathfrak R(F)$ that the decomposition  $\P(\mathfrak S)$ consists of  $m$ parts $G_{1,2}$, $G_{2,3}$,~\ldots, $G_{m,1}$ and $\R(G_{i,i+1})=Q_{i,i+1}$.  

We call the vertex~$p$  the {\it center}, and the vertices $q_1,\ldots,q_m$~--- {\it petals} of this flower. 
The set  $V(F)=\{p,q_1,\ldots,q_m\}$ is called the  {\it vertex set } 
of~$F$. All sets $Q_{i,j}=\{q_i, q_j, p\}$ are called {\it sets of a flower}~$F$.

We say, that the set $\mathfrak S$ {\it generates} the flower~$F$.
\end{defin}

The notations introduced above  will be standard for a flower. We will always write petals of a flower in cyclic order and consider their indexes as residues modulo number of petals. 
Set the notation $G_{i,j}=\cup_{x=i}^{j-1}\: G_{x,x+1}$ (index $x$ run over all residues from~$i$ to~$j-1$ in 
the cyclic order). We suppose, that $G_{x,x}=\varnothing$.

It is proved in~\cite[lemma~9]{k05}, that  the dependence graph of any set of  $k$-cutsets which generates a flower is connected. It is  also proved in~\cite{k05} (see theorem~6 and corollaries of it), that any set~$Q_{i,j}$
separates $G_{i,j}$ from $G_{j,i}$, and, moreover, $\P(Q_{i,j})=\{G_{i,j},G_{j,i}\}$ 
if $j\not\in \{i,i+1,i-1\}$.

We call the sets $Q_{1,2},Q_{2,3},\ldots,Q_{m,1}$  {\it boundaries}, and other sets  $Q_{i,j}$~--- {\it inner sets} of the flower~$F$.

Let $\P(F)= \{G_{1,2},G_{2,3},\ldots,G_{m,1}\}$ be the   {\it decomposition} of the graph $G$ by the flower~$F$. Obviously, no one of these parts is small. 

If a part $G_{i,i+1}$ is empty, then  $Q_{i,i+1}$ is not a cutset.
If  $G_{i,i+1}$ is nonempty, then $Q_{i,i+1}$ is a cutset,
 $G_{i+1,i}\in \P(Q_{i,i+1})$ and~$G_{i,i+1}$ is the union of all different from~$G_{i+1,i}$ parts of~$\P(Q_{i,i+1})$.

\begin{lem}
\label{lrr1}
Let a set  $\mathfrak S$ generate a flower $F$. Then the intersection of all cutsets of $\mathfrak S$ 
consists of one vertex~--- the center of the flower $F$.

\begin{proof} Let $F=\{p;q_1,\dots,q_m\}, \quad \mathfrak S=\{S_1,\dots,S_n\}, \quad P=\cap_{i=1}^n S_i.$
It is obvious from $\mathfrak S\subset \mathfrak R(F)$, that $P\ni p$. If $P\ne\{p\}$, then $P$ contains a petal of $F$, let $P\ni q_i$. Then all parts of $\P(\mathfrak S)=\P(F)$ contain  $q_i$. However, the set $G_{i+1,i-1} \not \ni q_i$ is a union of several parts of $\P(F)$. This contradiction finishes the proof.
\end{proof}
\end{lem}

\begin{rem}
\label{r30}
1) If a part~$G_{i,i+1}$ is empty, then according to the structure of a flower, described  above,  by corollary~\ref{l1c3} the vertices~$q_i$ and~$q_{i+1}$ are adjacent.

2) If both parts~$G_{i-1,i}$ and~$G_{i,i+1}$ are empty, then the vertex~$q_i$ is adjacent in the graph~$G$ to
the vertices $p,q_{i-1},q_{i+1}$ and only them.

3) If sets  $S=\{a,u,v\}$ and $T=\{a,x,y\}\in\mathfrak R_3(G)$ are dependent, then with the help of lemma~\ref{lds0} and corollary~\ref{l1c3} it is easy to see, that these two sets generate a flower with the center~$a$ and cyclic ordered petals $u,x,v,y$.
\end{rem}

A flower~$F$ can be generated by different sets of 3-cutsets, however, it is proved in~\cite[theorem~7]{k05},
that all such sets split the graph  identically~--- into the parts of~$\P(F)$. 
Moreover, it is proved in the mentioned theorem, that if  sets 
$\mathfrak S,\:\mathfrak T\in \mathfrak R_k(G)$ generate  flowers with the same center and the same set of petals, 
then these two flowers coincide  (i. e. the cyclic orders of petals in these flowers coincide).

\begin{figure}
\centerline{\includegraphics{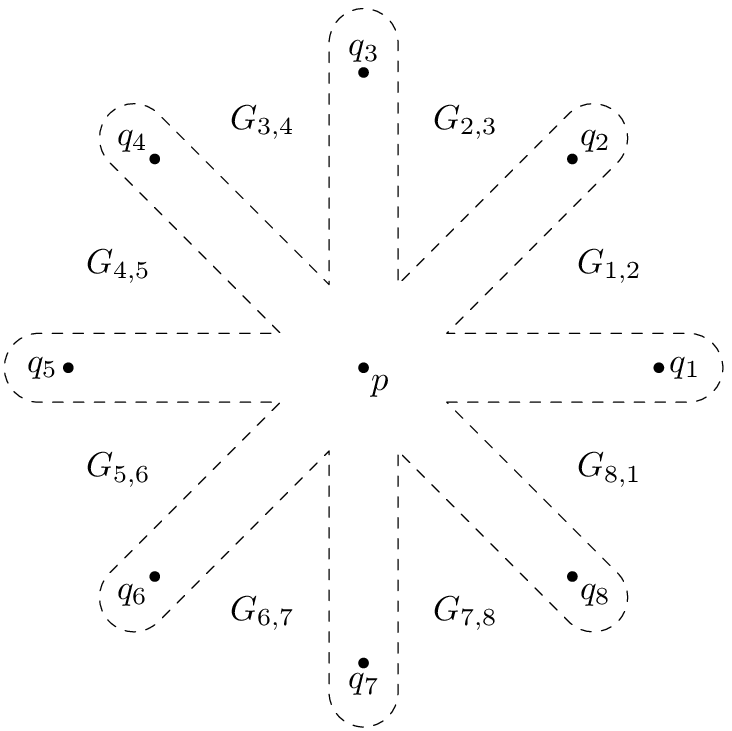}}
\caption{Decomposition of a graph by a flower with eight petals}\label{ris.1}
\end{figure}

On the figure~\ref{ris.1}  decomposition of a graph by a flower with eight petals it is shown.

In this sense a flower is alike a wheel, from which, according to the work~\cite{T2}, were ``originated''  all triconnected graphs.  

We need the following theorem, also proved in the  work~\cite{k05}.

\begin{thm}[{\cite[theorem~8]{k05}}]
\label{tr5}
For any set $\mathfrak S\subset \mathfrak R_3(G)$ two following statements are equivalent.

$1^\circ$ Every part of~$\P(\mathfrak S)$ contains at least three vertices.

$2^\circ$ Every dependence component of the set~$\mathfrak S$ either consists of one cutset, 
or generates a flower.
\end{thm}

\begin{cor}
\label{tr5c1}
If dependence graph of a set $\{S_1,S_2,\ldots,S_n\}\subset \mathfrak R_3(G)$ is connected and
 $\cap_{i=1}^n S_i \ne \varnothing$, then this set generates a flower.

\begin{proof}
Obviously, we can enumerate the cutsets of our set such, that for every $\ell\in\{1,\dots,n\}$ the dependence graph of the set  $\mathfrak S_\ell=\{S_1,S_2,\ldots,S_\ell\}$ is connected. Prove by induction on $\ell$, that the set  $\mathfrak S_\ell$ generates a flower.  The base for  $\ell=2$ is obvious by remark~\ref{r30}.

Induction step from~$\ell$ to~$\ell+1$.  Let the set $\mathfrak S_\ell$ generate a flower~$F=(p;q_1,\dots,q_m)$.  If  $S_{\ell+1}$ does not split any part of $\P(F)$, then $\P(\mathfrak S_{\ell+1})=\P(F)$. 
This decomposition contains no small parts and by theorem~\ref{tr5}  the step is proved.

Let   $S_{\ell+1}$ split some part of $\P(F)$.  It follows from remark~\ref{r30}, that the cutset~$S_{\ell+1}$ cannot split an empty  part of $\P(F)$. Let~$S_{\ell+1}$ split a nonempty part~$G_{i,i+1}$.  Then~$S_{\ell+1}$ 
is dependent with~$Q_{i,i+1}$. 
By lemma~\ref{lrr1} we have  $\cap_{i=1}^\ell S_i =\{p\}$, hence, $S_{\ell+1}\ni p$. Now it follows from  dependence of the sets $S_{\ell+1}$ and $Q_{i,i+1}$, that $S_{\ell+1}\cap Q_{i,i+1}=\{p\}$, the intersection $S_{\ell+1}\cap \I(G_{i,i+1})$ consists of a single vertex ~$x$, and vertices~$q_i$ and~$q_{i+1}$ lie in different parts of $\P(S_{\ell+1})$. 
Then by lemma~\ref{lds0} the set~$S_{\ell+1}$ splits the part~$G_{i,i+1}$ into two parts with boundaries
$\{q_i,p,x\}$ and $\{q_{i+1},p,x\}$. By corollary~\ref{l1c3}, both these parts are not small. Thus, $\P(\mathfrak S_{\ell+1})$ contains no small parts and by theorem~\ref{tr5} this set generates a flower. The 
induction step is proved. 
\end{proof}
\end{cor}

\begin{defin}
We say that a flower~$F$ {\it contains} a flower~$F'$, if they have common center  and
$V(F')\subset V(F)$.
We call a flower~$F$ {\it maximal}, if it is not contained in another flower.
\end{defin}

\begin{lem}
\label{l300}
Let $F=(p;q_1,\dots,q_m)$ be a maximal flower. Then the following statements hold.

$1)$ There is no set $T\in \mathfrak R_3(G) \setminus \mathfrak R(F)$, which contains~$p$ and is dependent
with at least one set of the flower~$F$.

$2)$ For each vertex  $v\in \I(G_{i,i+1})$ there exists a path between~$q_i$ and~$q_{i+1}$, which does not pass through~$v$, and all inner vertices of this path lie in $\I(G_{i,i+1})$.

\begin{proof} 
1) Assume the contrary. Then the set~$S$, dependent with a  set of the flower~$F$, must be dependent
with some cutset of~$\mathfrak R(F)$.  Hence, the dependence graph of the set $\mathfrak S=\mathfrak R(F) \cup \{S\}$ is connected. Moreover, every cutset of the set $\mathfrak S$ contains a vertex~$p$.
Then by corollary~\ref{tr5c1} there exists such a flower~$F'$ that
$\mathfrak R(F') \supset \mathfrak S$. Obviously, $F'$ contains $F$ and these flowers are different.
A contradiction with maximality of~$F$.

2)  Assume the contrary. Then it is easy to see, that the set $T=\{v,p,q_{i+2}\}$ separates~$q_i$ from~$q_{i+1}$, i.e., $T$ is a cutset dependent with~$Q_{i,i+1}$. 
A contradiction with item~1.
\end{proof}
\end{lem}

\begin{rem}
1) It is easy to reconstruct the center and petals of a flower $F$ by the set~$\mathfrak R(F)$.

2) Obviously, a flower~$F$ contains a flower~$F'$ if and only if $\mathfrak R(F)\subset \mathfrak R(F')$.

3) It is easy to derive from item~1 of lemma~\ref{l300} that every flower is contained in unique maximal
flower.
\end{rem}

\subsection{Vertex-edge cuts}

\begin{defin}
$1)$
Let a  {\it cutting set } be  any set  $M\in V(G)\cup E(G)$, for which  the graph $G-M$ is disconnected.
Obviously, every cutting set  of a triconnected graph contains at least three elements.

Denote by~$\mathfrak M_i(G)$ (where $i\in \{0,1,2,3\}$) the set, consisting of all cutting sets with  $i$ edges and  $3-i$ vertices of the graph~$G$.
Let $\mathfrak M(G)=\cup_{i=1}^3\: \mathfrak M_i(G)$, 
$\mathfrak M^+(G)=\mathfrak M(G)\cup\mathfrak M_0(G)$.
Note, that  $\mathfrak M_0(G)=\mathfrak R_3(G)$. 

All cutting sets of ${\mathfrak M}(G)$ we call {\it vertex-edge cuts}, or  simply {\it cuts}.

$2)$ Let $M, N\in \mathfrak M^+(G)$. If~$N$ contains all vertices of~$M$ and for every edge $e\in M$ the set $N$ contains either~$e$, or an end of $e$, we say that~$M$ {\it contains}~$N$ (or~$N$ {\it is contained} in~$M$).

If a cut  $M\in \mathfrak M(G)$ is not contained in any other cut of~$\mathfrak M(G)$, we call it
a {\it maximal} cut. 

We say that a cutting set   $M\in\mathfrak M^+(G)$ can be  {\it complemented} by an edge~$ab$ (or an edge~$ab$
{\it complements}~$M$), 
if  an end of $ab$ (let it be $a$) belongs to $M$ and after changing in $M$ the vertex $a$ to the edge~$ab$ we obtain a cut from $\mathfrak M(G)$.
\end{defin}

Note, that a cut is maximal if and only if it cannot be complemented by an edge.

\begin{rem}
\label{nr1} 
Let $M\in\mathfrak M(G)$.

$1)$ There is no vertex of $M$, which is incident to an edge of $M$.

$2)$ If two edges of the cut  $M$ have common end~$x$, then~$\{x\}$ is a connected component of the graph $G-M$. Indeed, otherwise after replacing these two edges by a vertex $x$ we obtain a cutting set of two elements, that is impossible in triconnected graph~$G$.

$3)$ We can consider the  cut  $M$ as a subgraph of~$G$ (containing vertices of $M$, edges of $M$ and ends of these edges). We denote by $V(M)$ the vertex set of this subgraph.
\end{rem}

Let a cut $M\in \mathfrak M(G)$ contain an edge~$x_1x_2$.
It is clear, that the graph $G-M$ has exactly two connected components: one of them contains~$x_1$, 
and the other contains~$x_2$.

\begin{defin}
Let  $x_1x_2\in M\in \mathfrak M(G)$ and the graph $G-M$ have two connected components~$H_1$ and~$H_2$, such that $x_1\in H_1$ and $x_2\in H_2$. Then we set the notations $G^M_i=V(G)\setminus H_{3-i}$ and $T^M_i=G^M_i\cap V(M)$ for $i\in\{1,2\}$.

We call the sets~$G^M_1$ and~$G^M_2$  {\it parts of $M$-decomposition} and use the notation
$\P(M)=\{G^M_1,G^M_2\}$. We call {\it interior} of the part~$G^M_i$  the set $\I(G^M_i)=G^M_i\setminus T^M_i$.
We call  {\it neighborhood} of the part~$G^M_i$  the set $\O(G^M_i)=G^M_i\cup V(M)$. Here $i\in\{1,2\}$.
We call {\it boundaries} of the cut $M$ the sets~$T^M_1$ and~$T^M_2$.

For any cut $M\in\mathfrak M(G)$ we shall use these notations. Every edge of the cut $M$ we shall write such that
it first end lies in~$G^M_1$, and second end lies in ~$G^M_2$.
\end{defin}

\begin{figure}
\centerline{\includegraphics{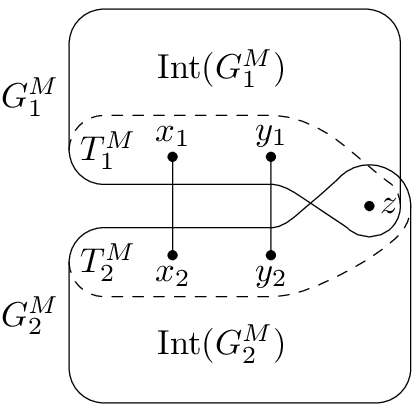}}
\caption{Decomposition of the graph by a cut from~$\mathfrak M_2(G)$}\label{ris.2}
\end{figure}

The sets~$T^M_i$, $G^M_i$ and $\I(G^M_i)$  are shown on fig.~\ref{ris.2} for a cut~$M\in\mathfrak M_2(G)$.

\begin{rem}
\label{rr5}
$1)$ Note, that the set~$G^M_i$ is obtained from the set~$H_i$ by adding of all vertices (but not edges!) of the cut~$M$. Thus, the part~$G^M_i$ contains all vertices of the cut~$M$ and exactly one end of each edge of the cut ~$M$. The set~$T^M_i$ also contains all vertices of the cut~$M$  
and exactly one end of each edge of the cut~$M$, and does not contain other vertices.
Hence, if   edges of~$M$ have no common ends, then $|T^M_1|=|T^M_2|=3$. If edges of~$M$ have common ends, then
by remark~\ref{nr1} one connected component of the graph $G-M$ (denote it by~$H_1$) consists of a single vertex. In this case $|T^M_1|=1$ and  $|T^M_2|=3$.

2) Note also, that the definition of a part of $M$-decomposition is compatible with the definition of a part of decomposition of the graph by a cutset:~$G^M_1$ and~$G^M_2$ are maximal (with respect to inclusion)
sets not splitted by~$M$.

3) Note, that a part of  $M$-decomposition, by contrast of a part of decomposition by a cutset, can consist of a single vertex. It happens when a cut~$M$ consists of three edges, incident to a vertex of degree~3.
\end{rem}

\begin{lem}
\label{ll31}
Let  $x\in M\in \mathfrak M^+(G)$, $xy \in E(G)$ and  
$H$ be a connected component of the graph $G-M$ containing~$y$.
Then the cutting set~$M$ can be complemented by an edge~$xy$ if and only if $y$ is the only
vertex of the component~$H$ which is adjacent to~$x$.

\begin{proof}
Let  $M'$ be a set obtained from~$M$ by replacing a vertex~$x$ by an edge~$xy$. 
If the graph $G-M'$ is disconnected, then it consists of exactly two connected components,
one of them contains the vertex~$x$, and the other component contains~$y$. Note, that all vertices 
of the component~$H$ lie in the same connected component of the graph~${G-M'}$. Thus, if the vertex~$x$
is adjacent in the graph~$G-M'$ to a vertex of~$H$, then the vertices~$x$ and~$y$ are connected in~${G-M'}$, i.e. the graph $G-M'$ is connected.

On the other side, if the vertex~$x$ is not adjacent to any vertex of the component~$H$, except~$y$, then there is no path between~$x$ and~$y$ in the graph $G-M'$, i.e. this graph is disconnected.
\end{proof}
\end{lem}

\begin{cor}
\label{ll31c1}
If $M\in\mathfrak M_2(G)$, then  there exists not more than one  edge which complement~$M$.

\begin{proof}
Let~$x$ be the only vertex of the cut~$M$ and~$H_1$, $H_2$ be connected components of the graph~$G-M$. 
Then by lemma~\ref{ll31} there is not more than one edge from $x$ to each of these components by which $M$ can be complemented. Obviously,  $V(G)=H_1\cup H_2\cup \{x\}$. Since the graph~$G$ is triconnected, $d(x)\ge3$. 
Thus the vertex~$x$ cannot be adjacent to exactly one vertex of the component~$H_1$ and exactly one
vertex of the component~$H_2$. Hence, the cut~$M$ can be complemented by not more than one edge.
\end{proof}
\end{cor}

\begin{cor}
\label{ll31c2}
Let cutsets $S,T\in\mathfrak R_3(G)$ be dependent and $\{x,y\}\in\P(\{S,T\})$.
Then each of the sets~$S$ and~$T$ can be complemented by the edge~$xy$.

\begin{proof}
By corollary~\ref{l1c3} we have $S \cap T = \varnothing$.
Without loss of generality we may suppose, that $x\in S$, $y\in T$. Let $x\in F\in\P(T)$, 
$y\in H\in\P(S)$. By corollary~\ref{l1c2} the vertices~$x$ and~$y$ are adjacent and $y$ is the only vertex
of the part $H$ adjacent to~$x$. Then by lemma~\ref{ll31} the set~$S$ can be complemented by the edge~$xy$. 
Similarly, the set~$T$ can be complemented by the edge~$xy$. 
\end{proof}
\end{cor}

\begin{defin}
We call  a cut $M\in\mathfrak M(G)$ {\it nondegenerate}, if 
$\I(G^M_1)\ne\varnothing$ and $\I(G^M_2)\ne\varnothing$. Otherwise, we call this cut  {\it degenerate}.

We call a cutting set $M\in\mathfrak M^+(G)$ {\it trivial}, if 
one   connected component of the graph $G-M$ consists of a single vertex, and {\it  nontrivial} otherwise.
\end{defin}

A case  $\I(G^M_1)=\I(G^M_2)=\varnothing$ is not interesting, because in this case the graph $G$ contains not more than 6 vertices. Further we suppose, that for every degenerate cut $M$ exactly one of the sets
 $\I(G^M_1)$ and $\I(G^M_2)$ is empty.

\begin{figure}
\centerline{\includegraphics{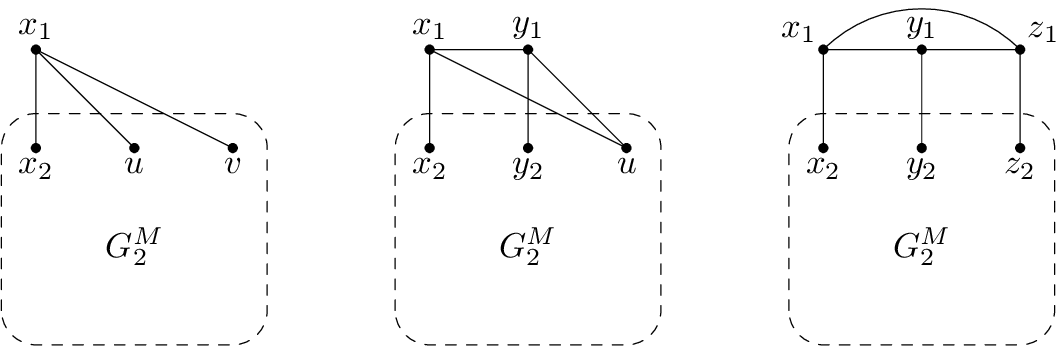}}
\caption{Degenerate cuts.}\label{ris.3}
\end{figure}

\begin{rem}
\label{nr2}
1) A degenerate cut containing exactly one edge is trivial.
Indeed, if a cut $M=\{u,v,x_1x_2\}$ is degenerate ($\I(G_1^M)=\varnothing$), then since the graph $G$ is triconnected, the vertex~$x_1$ is adjacent to the vertices  $x_2,u,v$ and only them. 
Thus the set  $T_2^M=\{x_2,u,v\}$ separates the vertex~$x_1$ from other vertices of the graph. Sets $\{x_1u,x_1v,x_1x_2\}$ and~$\{u,x_1v,x_1x_2\}$ also are  cuts. By remark~\ref{nr1}, if a cut contains 
adjacent edges, it is trivial.

2) The structure of a degenerate nontrivial cut  also can be simply described. 
 It follows from written above, that such cut contains  more than one edge.

If a cut $M=\{u,x_1x_2,y_1y_2\}$ is degenerate ($\I(G_1^M)=\varnothing$) and nontrivial, then the vertex~$x_1$ is adjacent to the vertices $y_1,x_2,u$ and only them, and the vertex~$y_1$ is adjacent to~$x_1,y_2,u$ and only them.

If a cut $M=\{x_1x_2,y_1y_2,z_1z_2\}$ is degenerate and nontrivial (again $\I(G_1^M)=\varnothing)$, then the vertices~$x_1$, $y_1$, $z_1$ are pairwise adjacent and except this edges the vertices~$x_1$, $y_1$, $z_1$  are incident to  edges of the cut $M$ and only them.

Degenerate cuts with one, two and three edges are shown on figure~\ref{ris.3}.
\end{rem}

\begin{lem}
\label{l30}
Let $M\in \mathfrak M(G)$ be a nontrivial cut. Then the following statements hold.

$1)$ Every set, which contains all vertices of $M$ and exactly one end of each edge of $M$ and differs 
from~$T^M_1$ and~$T^M_2$, is a cutset. 
Moreover, this cutset splits the graph into two parts, one of which contains~$G^M_1$, and the other contains~$G^M_2$.

$2)$ If $\I(G_2^M)\ne\varnothing$, then $T^M_2$ is a cutset. Moreover,  $\O(G^M_1)\in \P(T^M_2)$ and~$G^M_2$ is a union of several parts of~$\P(T^M_2)$. If the cut~$M$ is nondegenerate, then both~$T^M_1$ and~$T^M_2$
are cutsets.

\begin{proof}
1) Let~$R$ be any such set. Since~$R$ does not coincide with~$T^M_1$ and~$T^M_2$, then the sets~$G^M_1\setminus R$ and~$G^M_2\setminus R$ are nonempty. Obviously,~$R$ separates these two sets from each  other, thus, $R$ is a cutset. 

Let us prove, that  vertices of the set $G^M_1\setminus R$ are connected in the graph~${G-R}$. Indeed,
let~$x_1\in T^M_1\setminus R$.  
Then the cut~$M$ contains an edge~$x_1x_2$ and $x_2 \in R$. Since the cut~$M$ is nontrivial, then~$x_1$  is the only vertex outside the part~$G^M_2$ which is adjacent to~$x_2$. However, each connected 
component of the graph~$G-R$ must contain a vertex adjacent to~$x_2$. Since any connected component of~$G-R$ which contains a vertex of $G^M_1\setminus R$ contain no vertices of~$G^M_2$, then all vertices of the set~${G^M_1\setminus R}$ are  in the same connected component of the graph~${G-R}$.

Analogously, all vertices of the set~${G^M_2\setminus R}$ are  in the same connected component of the graph~${G-R}$. Hence there are exactly two connected components in the graph~${G-R}$.

2) It is easy to see, that the set~$T_2^M$ separates~$G_2^M$ from~$\O(G_1^M)$. Hence,~$T_2^M$ is a cutset.
The fact, that  vertices of the set   $G^M_1\setminus T_2^M$ are  connected in the graph $G-R$, is proved as well as in item~1. It is clear, that a part of $\P(T_2^M)$ 
containing~${G^M_1\setminus T_2^M}$ is~$\O(G^M_1)$.
\end{proof}
\end{lem}

\begin{cor}
\label{l30c1}
Let $M\in\mathfrak M(G)$, a cutset $S\in\mathfrak R_3(G)$ be dependent with $T^M_2$, $S\cap G^M_2=\{x\}$ 
and~$S$ separates a vertex~$y\in T^M_2$ from other vertices of~$T^M_2$. Then the cut~$M$ can be complemented by the edge~$xy$.

\begin{proof}
Let $x\in H\in\P(T^M_2)$ (it is easy to check, that~$H=G^M_2$) and $y\in F\in\P(S)$.
By condition, $S\cap H=\{x\}$ and $T^M_2\cap F=\{y\}$.
Then by corollary~\ref{l1c2}, we have that $\{x,y\}\in\P(\{S,T^M_2\})$, the vertices~$x$ and~$y$ are adjacent
and $x$ is the only vertex of the part~$G^M_2$ adjacent to~$y$. Hence, by lemma~\ref{ll31} the edge~$xy$
complements the cut~$M$.
\end{proof}
\end{cor}

\begin{rem}
\label{nr3}
Obviously, all cutsets of lemma~\ref{l30} are contained in the cut~$M$.
Moreover, all these sets, except~$T^M_1$ and~$T^M_2$, are pairwise dependent. The cutsets~$T^M_1$ and~$T^M_2$ are independent with each other and with all other cutsets, contained in~$M$.
\end{rem}

\begin{defin}
The cutsets described in item~1 of lemma~\ref{l30} we call 
{\it inner sets} of the cut~$M$. The set consisting of all inner sets of~$M$ we denote by~$\mathfrak R(M)$.
\end{defin}

\subsubsection{Singular edges}

\begin{lem}
\label{ll32}
Let cuts $M_1,M_2\in\mathfrak M_2(G)$ have two common edges. Then there is a cut $M\in\mathfrak M_3(G)$ containing both~$M_1$ and~$M_2$.

\begin{proof}
Let $M_1=\{x_1x_2,y_1y_2,z_1\}$ and~$z_2$ be the only vertex of~$M_2$. 
Without loss of generality we may suppose, that $z_2 \in G^{M_1}_2$. Then~$M_2$
does not split~$G^{M_1}_1$. In particular, $M_2$ does not split~$\{x_1,y_1,z_1\}$. Thus we may suppose, that $M_2=\{x_1x_2,y_1y_2,z_2\}$ and $z_1 \in G^{M_2}_1$. 

Note, that then the cutset~$T^{M_1}_2=\{x_2,y_2,z_1\}$ separates~$z_2$ from  $\{x_1,y_1\}$, 
and the cutset~$T^{M_2}_1=\{x_1,y_1,z_2\}$ separates~$z_1$ from $\{x_2,y_2\}$. Hence, the cutsets~$T^{M_1}_2$ 
and~$T^{M_2}_1$ are dependent and by corollary~\ref{l30c1} we can complement the cut~$M_1$ by the edge~$z_1z_2$ and obtain desired cut $M=\{x_1x_2,y_1y_2,z_1z_2\}$.
\end{proof}
\end{lem}

\begin{cor}
\label{ll32c1}
Two maximal cuts cannot  have  more than one common edge.

\begin{proof}
Let cuts $M_1,M_2\in\mathfrak M(G)$ have two common edges. If $M_i\in\mathfrak M_2(G)$, then we set $M_i'=M_i$, 
else (when $M_i\in\mathfrak M_3(G)$), we obtain  $M_i'$ from $M_i$ 
replacing the edge which is not in $M_{3-i}$ by one of its ends.
Clearly, we can perform this replacement such that cuts $M_1'$ and~$M_2'$ would be different.
Then by lemma~\ref{ll32} there is a cut $M\in\mathfrak M_3(G)$, containing both cuts  $M_1'$ and~$M_2'$. 
However, by corollary~\ref{ll31c1} each of the cuts~$M_1'$ and~$M_2'$ can be contained 
in not more than one cut from~$\mathfrak M_3(G)$. Hence, each of the cuts $M_1$ and~$M_2$ is contained in $M$ or coincide with~$M$, i.e. at least one of the cuts $M_1$ and~$M_2$ is not maximal.
\end{proof}
\end{cor}

\begin{defin}
We call an edge  $e\in E(G)$ {\it singular}, if there exist different vertices $u,v,t,w\in V(G)$ such that $\{u,v,e\},\{t,w,e\}\in \mathfrak M(G)$.
\end{defin}

\begin{rem}
\label{r22}
Let $\{a_1a_2,b_1b_2,c_1c_2\}\in\mathfrak M_3(G)$.
It is easy to see, that all edges~$a_1a_2$, $b_1b_2$, $c_1c_2$ are singular.
\end{rem}

\begin{defin}
Let  $M=\{u,v,x_1x_2\},\: N=\{t,w,x_1x_2\}\in\mathfrak M_1(G)$, and the vertices $u,v,t,w$ are different. 
We call the cuts~$M$ and~$N$  {\it independent}, if one of the parts~$G^M_1$ and~$G^N_1$ contains the other. Otherwise, we call these cuts  {\it dependent}.
\end{defin}

Clearly, the cuts~$M$ and~$N$ are independent if and only if either $t,w\in \I(G^M_1)$, or $t,w\in \I(G^M_2)$.

\begin{lem}
\label{ll33}
Let cuts $M,N\in\mathfrak M_1(G)$ be dependent and both contain an edge~$x_1x_2$. 
Then there exists a cut $\{x_1x_2,y_1y_2,z_1z_2\}\in\mathfrak M_3(G)$ such that
$M=\{x_1x_2,y_1,z_2\}$, $N=\{x_1x_2,y_2,z_1\}$.

\begin{proof}
Let $M\cap G^N_1 = \{y_1\}$, $M\cap G^N_2 = \{z_2\}$, $N\cap G^M_1 = \{z_1\}$,
$N\cap G^M_2 = \{y_2\}$. Consider  sets  $T^M_2=\{x_2,y_1,z_2\}$ and~$T^N_1=\{x_1,y_2,z_1\}$. 
Note, that  $T^M_2$ separates the vertex~$y_2$ from~$\{x_1,z_1\}$, and $T^N_1$ separates the vertex~$y_1$ from $\{x_2,z_2\}$.
Thus, by corollary~\ref{l30c1} both cuts~$M$ and~$N$ can be complemented by the edge~$y_1y_2$.
Obtained cuts $\{x_1x_2,y_1y_2,z_1\}$ and $\{x_1x_2,y_1y_2,z_2\}$ by lemma~\ref{ll32}
can be complemented by the edge~$\{z_1z_2\}$. As a result, we obtain desired cut~$\{x_1x_2,y_1y_2,z_1z_2\}$.
\end{proof}
\end{lem}

\begin{thm}
\label{t30}
For vertices  $x_1,x_2\in V(G)$ two following statements are equivalent.

$1^\circ$ Vertices~$x_1$ and~$x_2$ are adjacent, $x_1x_2$ is a singular edge.

$2^\circ$ There exist dependent cutsets $S,T\in\mathfrak R_3(G)$ such that
 $x_1\in S$, $x_2\in T$ and $\{x_1,x_2\}\in\P(\{S,T\})$.

\begin{proof}
$2^\circ \Rightarrow 1^\circ$. Let the condition~$2^\circ$ hold. Then, by corollary~\ref{ll31c2}
both cutsets~$S$ and~$T$ can be complemented by the edge~$x_1x_2$. Hence, the edge~$x_1x_2$ is singular.

$1^\circ \Rightarrow 2^\circ$. If the  condition~$1^\circ$ holds, then there exist different vertices~$u,v,t,w\in V(G)$ such 
that $M=\{u,v,x_1x_2\}, N=\{t,w,x_1x_2\}\in \mathfrak M(G)$.
Consider two cases.

\q{a}. Suppose, that~{\it $M$ and~$N$ are independent}. Without loss of generality
we may assume, that $G^M_1\supset G^N_1$ and $G^M_2\subset G^N_2$. 
Consider disjoint sets~$T^M_1$ and~$T^N_2$. In our case $t,w\in \I(G^M_1)$, $x_2\not\in G^M_1$.
Then by lemma~\ref{l30} the set~$T^M_1$ is a cutset and separates~$x_2$ from~$\{t,w\}$.
Similarly,~$T^N_2$ is a cutset and separates~$x_1$ from~$\{u,v\}$.
Hence, the cutsets  $T^M_1=\{x_1,u,v\}$ and $T^N_2=\{x_2,t,w\}$ are dependent and by corollary~\ref{l1c2}
we have $\{x_1,x_2\}\in\P(\{T^M_1,T^N_2\})$.

\q{b}. Suppose, that~{\it $M$ and~$N$ are dependent}. 
Then by lemma~\ref{ll33} there exists a cut $\{x_1x_2,y_1y_2,z_1z_2\}\in\mathfrak M_3(G)$ such that 
 $M=\{x_1x_2,y_1,z_2\}$, $N=\{x_1x_2,y_2,z_1\}$. Consider sets~$S=\{x_1,y_2,z_2\}$ and $T=\{x_2,y_1,z_1\}$. 
By lemma~\ref{l30} we have $S,T\in\mathfrak R_3(G)$, moreover,~$S$ separates the vertex~$x_2$ 
from~$y_1,z_1$, and~$T$ separates the vertex~$x_1$ from~$\{y_2,z_2\}$. Thus, the sets~$S$ and~$T$ are dependent and by corollary~\ref{l1c2} we have $\{x_1,x_2\}\in\P(\{S,T\})$.
\end{proof}
\end{thm}

\subsection{Trivial cutsets and triple cuts}

Remind, that a cut~$M$ is called  {\it trivial}, if one connected component of the graph~${G-M}$ 
consists of one vertex. Obviously, degree of this vertex is~3. In this section we study trivial cutsets.

\begin{defin}
Let $T\in \mathfrak R_3(G)$ be a trivial cutset, separating a vertex~$a$ from other vertices.
Let a set $S\in\mathfrak R_3(G)$, containing the vertex $a$, be such that~$|\P(S)|=3$ and interior
of every part of~$\P(S)$ contains a vertex of the set~$T$. 
Then we say, that the trivial set~$T$ is {\it subordinated} to the set~$S$. 
Denote by~$\mathfrak D$ the set of all  3-cutsets, which have a subordinated  trivial cutset.
\end{defin}

\begin{rem}
It is clear by definition, that if $S\in \mathfrak D$, then there exists a vertex $a\in S$ of degree 3.
\end{rem}

\begin{lem}
\label{l3v0}
$1)$ If a cutset $S\in \mathfrak R_3(G)$ splits $G$ into more than three parts, then $S$ is independent with all cutsets of~$\mathfrak R_3(G)$.
If a cutset $S\in \mathfrak R_3(G)$ splits $G$ into three parts and $S$ is dependent with 
$T\in \mathfrak R_3(G)$, then the cutset~$T$ is trivial and subordinated to~$S$.

$2)$ A trivial cutset can be subordinated to not more than one cutset.

\begin{proof} 
1) Let cutsets $S,T\in\mathfrak R_3(G)$ be dependent. Then by corollary~\ref{l1c1} we have 
$|\P(S)|\le3$ and $|\P(T)|\le3$.
Let $\P(S)=\{H_1,H_2,H_3\}$. If  $|\P(T)|=3$, then, by remark~\ref{r0} we have 
$|V(G)|=6$, this case is not interesting for us.
It is enough to consider the case $|\P(T)|=2$. Let $\P(T)=\{F_1,F_2\}$. 
We enumerate the parts such that $\I(F_1)\cap S= \{a\}$, $|\I(F_2)\cap S|=2$. 
Then by corollary~\ref{l1c2} the set~$S$ splits~$F_1$ into three empty parts.
Hence, $\I(F_1)= \{a\}$, i.e.~$T$ is a trivial cutset subordinated to~$S$.

2) Let~$T$ be a trivial cutset subordinated to both cutsets~$S$ and~$S'$. 
By item 1 then $|\P(S)|=|\P(S')|=3$, hence, the cutsets~$S$ and~$S'$ are independent. 
Let $A\in \P(S)$ be a part containing~$S'$.
Two vertices of the cutset~$T$ lie outside~$A$. Clearly,~$S'$ cannot separate these two vertices from each other. We have a contradiction.
\end{proof}
\end{lem}

Further in this section we consider a cutset $S\in\mathfrak D$.

\begin{defin}
Consider a vertex $a\in S$ of degree~3. Three vertices adjacent to~$a$ form a cutset, separating $a$ from other vertices.  We denote this cutset by~$T(a)$ and call it {\it neighborhood} of the vertex~$a$.
\end{defin}

It easy to see, that if $a\in S$ and $d(a)=3$, then $T(a)$ is a trivial cutset, subordinated to~$S$. 
Moreover, the set~$S$ can be complemented by any edge~$aa_i$ where~$a_i\in T(a)$.

Let $\P(S)=\{A_1,A_2,A_3\}$. 
We replace in $S$ every vertex~$a$ of degree~3 by an edge connecting $a$ with the vertex of 
$T(a)\cap \I(A_i)$ and denote obtained cut by~$M_i$.
It follows from above, that  $M_1,M_2,M_3\in\mathfrak M(G)$. 

The cuts $M_1,M_2,M_3$ may be not maximal. If $M_i$ is contained in a cut from~$\mathfrak M_3(G)$, denote this cut by $M_i'$ (obviously, this cut is unique).
In all other cases (among them the case when $M_i\in\mathfrak M_1(G)$ is contained in a cut from $\mathfrak M_2(G)$) we set $M_i'=M_i$.

Obviously,  $V(M_i)\subset V(M_i')\subset A_i$. Moreover, the set~$S$ is a bound of both cuts~$M_i$ 
and~$M_i'$. By lemma~\ref{l30} we have, that $A_{i+1}\cup A_{i+2}\in\P(M_i)$ and 
$A_{i+1}\cup A_{i+2}\in\P(M'_i)$ (the numeration is cyclic modulo 3). 
Denote the {\it other part} of $\P(M_i)$ by~$B_i$, and its boundary by~$T_i$. We denote by~$B_i'$ the part 
of~$\P(M_i')$, contained in~$B_i$, and its boundary denote by~$T_i'$. 
The neighborhood of $B_i$ as a part of  $\P(M_i)$ and the  neighborhood of $B'_i$ as a part of  $\P(M'_i)$
are defined. It is easy to see, that $\O(B_i)=\O(B_i')=A_i$. 

Note, that the cut~$M_i$ (and, consequently, the cut~$M_i'$) can be trivial. 
In this case~$|B_i'|=1$. Moreover, if in this case all vertices of the cutset~$S$ are of degree~3, then 
also~$|B_i|=1$.

\begin{defin}
We call  $F=M_1\cup M_2\cup M_3$ a {\it triple cut}, and 
$\O(F)=V(M_1')\cup V(M_2')\cup V(M_3')$~--- its  {\it neighborhood}.
The set~$S$ we call a {\it line} of triple cut. All inner cutsets of the cuts~$M_i$, the set~$S$ and all
cutsets subordinated to $S$ we call {\it inner sets} of this triple cut. We call the sets  $T_1$, $T_2$, $T_3$  {\it boundaries} of our triple cut and the sets  $T'_1$, $T'_2$, $T'_3$ ---  boundaries of its neighborhood.
Set $V(F)=V(M_1)\cup V(M_2)\cup V(M_3)$ and $\P(F)=\{B_1,B_2,B_3\}$.
\end{defin}

\begin{figure}
\centerline{\includegraphics{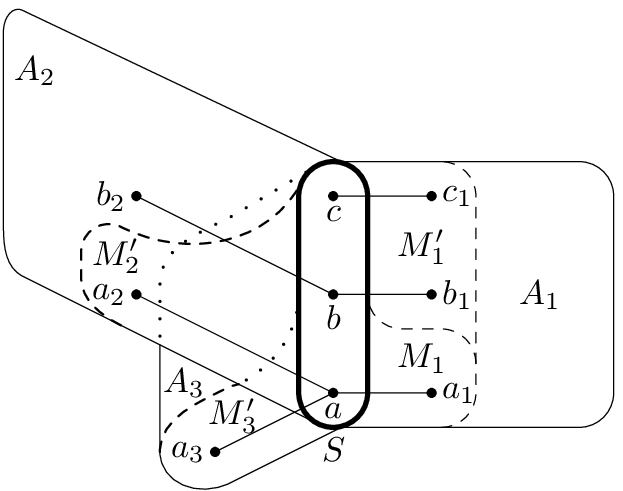}}
\caption{A triple cut with one trivial cutset}\label{ris.4}
\end{figure}

An example of a triple cut is shown on figure~\ref{ris.4}. In this example the line  $S=\{a,b,c\}$ of this 
triple cut has one subordinated cutset $T(a)=\{a_1,a_2,a_3\}$.
Here $M_1=\{aa_1,b,c\}$, $M_1'=\{aa_1,bb_1,cc_1\}$, $M_2=M_2'=\{aa_2,b,c\}$, $M_3=M_3'=\{aa_3,b,c\}$. Note, that the cut  $M_2'$ coincides with $M_2$, in spite of it is  not maximal: the cut $M_2$ can be complemented by the 
edge~$bb_2$.

\begin{rem}
A nontrivial definition of the cut~$M_i'$ is concerned with our aim.
On one side, we want to give  the most simple description of the parts
of decomposition of the graph by a set of cutsets, contained in the neighborhood of a triple cut, which are 
not its boundaries. On the other side we want each 3-cutset not contained in~$\O(F)$ to separate from the neighborhood not more than one vertex. Later we shall show, that
our definition of the neighborhood satisfies both these conditions.
\end{rem}

\section{Further properties of basic structures}

The structures described in previous section are generated by sets of dependent cutsets. In this section we describe their connection with each other and with other cutsets.

\subsection{Inner cutsets}

\begin{lem}
\label{lor12}
Let~$S$ be a set of three petals of a flower~$F$. Then~$S$ is not a cutset.  

\begin{proof}
Let $\I(G_{i,i+1})\ne\varnothing$. Then all vertices of the set $G_{i,i+1}\setminus S$
are connected in~$G-S$, and~$p$ is among them. Thus all vertices of~$G-S$ are connected, except, may be,
some petals of the flower, not belonging  to nonempty parts of~$\P(F)$. 
But any such petal~$q_j$ belongs to two empty parts of~$\P(F)$ and, by remark~\ref{r30}, is adjacent to~$p$. Thus, the graph~$G-S$ is connected.
\end{proof}
\end{lem}

\begin{cor}
\label{lor12c1}
Any $3$-cutset,  which is a subset of the vertex set of a flower, contains its center,
i.e. it is either an inner set of this flower, or its boundary.
\end{cor}

\begin{lem}
\label{l31}
Let $M\in \mathfrak M(G)$ be a nontrivial cut, containing an edge~$x_1x_2$. 
Then there is no cutset $S\in \mathfrak R_3(G)$, which contains both~$x_1$ and~$x_2$.

\begin{proof}
Suppose the converse, let $x_1,x_2\in S\in\mathfrak R_3(G)$.
Without loss of generality suppose, that $S=\{x_1,x_2,t\}$, where $t\in G^M_1$.
Since $\I(G^M_2)$ is a union of interiors of several parts of~$\P(T^M_2)$ and~$S\cap \I(T^M_2)=\varnothing$, 
then by lemma~\ref{lds1} all vertices of the set~$G^M_2\setminus S$  are  connected in the graph~$G-S$. 
Moreover, they are also connected with vertices of the set~$T^M_1\setminus S$, because each vertex of this set either belongs to the set $T^M_2\setminus S$, or is adjacent to a vertex of this set.

Consider any vertex~$w\in \I(G^M_1)\setminus S$, if this set is nonempty.
Since~$\I(G^M_1)$ is a union of interiors of several parts of~$\P(T^M_1)$, by Menger's theorem there
exist three vertex-disjoint paths inside~$G^M_1$ from $w$ to three vertices of the set~$T^M_1$. 
Since~$|S\cap G^M_1|=2$, at least one of these paths omit~$S$. Hence the vertex $w$ is connected in the 
graph~$G-S$  with the set~$G^M_2\setminus S$. That means the graph~$G-S$ is connected, we have a contradiction.
\end{proof}
\end{lem}

\begin{cor}
\label{l31c1}
Any $3$-cutset, which is contained in the vertex set of a cut, contains all its vertices and exactly one vertex of each edge of this cut. Hence, this cutset is an inner set of this cut or its boundary.
\end{cor}

\begin{lem}
\label{ll11}
For a triple cut $F=M_1\cup M_2\cup M_3$ the following statements hold.

$1)$ Any $3$-cutset, which is contained in $V(F)$ is an inner set of this triple cut or its boundary.

$2)$ Any $3$-cutset which is contained in the neighborhood of a triple cut $F$ either is subordinated to the line of $F$, or is contained in one of the cuts $M'_1$, $M'_2$, $M'_3$.

\begin{proof}
1) A cutset dependent with the line of a triple cut~$F$ by lemma~\ref{l3v0} is subordinated to it. Hence, 
this cutset is an inner set of~$F$. The cutset, independent with the line of~$F$ is contained in one of the sets~$V(M_i)$
(where $i\in\{1,2,3\}$). Then by corollary~\ref{l31c1} this cutset is either an inner set or a boundary of the cut~$M_i$, and by definition it is an inner set or a boundary of the triple cut~$F$.

2) Similarly to item 1, a set independent with the line of~$F$ is contained in one of the 
sets~$V(M_i')$ and, consequently, is contained in the cut~$M_i'$.
\end{proof}
\end{lem}

\subsection{Connection between flowers and cuts}

In this section we consider a question, in what cases vertex sets of a cut and a flower  coincide, or one of them is a subset of the other.

\begin{defin}
We  say that a flower~$F$ {\it is contained} in a cut~$M$, if $V(F)\subset V(M)$.
We  say that a cut~$M$ {\it is contained} in a flower~$F$, if $V(M)\subset V(F)$.
\end{defin}

\begin{lem}
\label{ll12}
A flower contained in a cut has exactly $4$ petals and two non-neighboring empty parts.

\begin{proof}
Let~$M$ be a cut and $F$ be a flower such that $V(F)\subset V(M)$.
Obviously, $|V(F)|\le |V(M)|\le 6$. Moreover, if $V(F)=6$, then $M\in\mathfrak M_3(G)$,
$V(F)=V(M)$ and cut~$M$ is nontrivial. Then the center of~$F$ is an end of an edge of the cut~$M$, and the other end of this edge is a petal of~$F$. Hence, both ends of this edge belong to a cutset of~$F$. That is impossible
by lemma~\ref{l31}.

Hence, $V(F)=5$ and the flower~$F$ has exactly  4 petals. Similarly to written above, the center of $F$ cannot 
be connected with its petal by an edge of the cut~$M$. Thus, petals of $F$ form two pairs of vertices, connected by edges of~$M$. Obviously, the petals, which are ends of an edge of $M$ are neighboring, because non-neighboring
petals are not adjacent. Let them be $q_1$ and~$q_2$. Then the  part $G_{1,2}$ is empty, since otherwise $q_1$ and~$q_2$ are connected by a path inside~$G_{1,2}$, and this path does not contain edges of $M$.
Clearly, this is impossible.
\end{proof}
\end{lem}

Therefore, only a flower with 4 petals can be contained in a cut. It is easy to see, that for every nontrivial
cut of~$\mathfrak M_2(G)$ its two inner sets generate a flower with 4 petals, which vertex set coincide with the vertex set of considered cut. This is the only case, when  vertex sets of a cut and of a flower coincide.
For every nontrivial cut~$M\in\mathfrak M_3(G)$ there are 6 cuts of~$\mathfrak M_2(G)$ contained in $M$ and 6 flowers correspondent to these cuts. These 6 flowers are contained in~$M$.

\begin{defin}
We call a flower~$F$ {\it nondegenerate}, if it is not contained in any cut~$M\in\mathfrak M_3(G)$, and {\it degenerate} otherwise.
\end{defin}

Now consider more often situation, when a flower contains  a cut.

\begin{lem}
\label{l32}
Let $F=(p;q_1,\dots,q_m)$ be a maximal flower and $\{q_i,p,q_jx\}\in \mathfrak M(G)$. 
Then one of the following two statements holds.

$1^\circ$ The vertex~$x$ is a petal of $F$, neighboring with~$q_j$, and $\{q_j,p,x\}$ is an empty part of~$\P(F)$.

$2^\circ$ The conditions $\{i,j\}=\{k,k+1\}$ and $x\in \I(G_{k,k+1})$ hold. Moreover, if $|\P(Q_{k,k+1})|=2$, 
then the vertices~$q_k$ and~$q_{k+1}$ are adjacent.

\begin{proof} 
Clearly, both ends of the edge~$q_jx$ lie in the same part of~$\P(F)$. Without loss of generality we may suppose
that $x\in G_{j,j+1}$. Let  $x\in H\in\P(Q_{j,i})$. Note, that if $|\P(Q_{j,i})|=2$, then 
$H=G_{j,i}$ (this condition for $i\ne j+1$ certainly holds). 
By lemma~\ref{ll31} we have, that~$q_j$ is not adjacent to different from~$x$ vertices of~$\I(H)$. 

Let $x\ne q_{j+1}$. If  $i\ne j+1$, then it is easy to see, that the set $T=\{q_i,p,x\}$ 
separates~$q_j$ from~$q_{j+1}$ and, consequently, $T$ is a cutset. By lemma~\ref{l300}, this contradicts
the maximality of the flower~$F$. Thus, $i=j+1$. Now note, that if $|\P(Q_{j,j+1})|=2$, then~$q_j$ 
is not adjacent to vertices of~$\I(G_{j,j+1})$ different from~$x$. 
If~$q_j$ and~$q_{j+1}$ are not adjacent, then, clearly, the set $T_1=\{q_{j+2},p,x\}$ separates~$q_j$ from~$q_{j+1}$ and, consequently, $T_1$ is a cutset. That also contradicts  maximality of the flower~$F$.
Hence, in this case the vertices~$q_j$ and~$q_{j+1}$ are adjacent.  

Let $x=q_{j+1}$. Then~$q_j$ is not adjacent to any vertex of~$\I(G_{j,j+1})$, whence it follows that the part~$G_{j,j+1}$ is empty.
\end{proof}
\end{lem}

\begin{cor}
\label{l32c1}
$1)$ Let~$M$ be a nontrivial cut and $F=(p;q_1,\ldots,q_m)$ be a flower such that $V(M)\subset V(F)$.
Then $p\in M$ and each edge of the cut~$M$ takes the form~$q_jq_{j+1}$, where  $G_{j,j+1}$ is an empty part.

$2)$ Let $F=(p;q_1,\dots,q_m)$ be a flower and $\I(G_{i,i+1})=\varnothing$.
Then $\{q_j,p,q_iq_{i+1}\}\in \mathfrak M_1(G)$ for every $j\not\in\{i,i+1\}$.
Moreover, if $\I(G_{j,j+1})=\varnothing$ and $i\ne j$, then 
$\{q_jq_{j+1},p,q_iq_{i+1}\}\in \mathfrak M_2(G)$.

\begin{proof} 
1) Since at least one boundary of the cut~$M$ by lemma~\ref{l30} is a cutset, by
lemma~\ref{lor12} this boundary contains the center of the flower~$F$. Consequently, $p\in V(M)$.
Suppose that $px\in M$. Clearly, there exist a cutset of flower~$F$ containing both~$p$ and $x$, that
contradicts lemma~\ref{l31}. Hence, $p\in M$.

Let $xy\in M$. Consider a maximal flower~$F'$ which contains~$F$. Since~$x$  and~$y$ are petals of~$F$, they are 
also petals of~$F'$. By lemma~\ref{l32}, the petals~$x$ and~$y$ must be neighboring in the flower~$F'$ 
(hence, in the flower~$F$ too), and correspondent to $x$ and~$y$ part of $\P(F)$ is empty.

2) It is easy to verify, that each of these sets separates~$q_i$ from~$q_{i+1}$.
\end{proof}
\end{cor}

\begin{rem}
A trivial cut can be contained in a flower, if this flower has two neighboring empty parts. Then edges connecting their common petal with neighboring petals and with the center of flower form a trivial cut. It is easy to check, that any trivial cut, contained in a flower is of this type.
\end{rem}

\begin{defin}
Let $G_{i,i+1}\in\P(F)$. Define a set $M_{i,i+1}$ as follows.

$1^\circ$~$p\in M_{i,i+1}$;
 
$2^\circ$~$q_{i-1}q_i\in M_{i,i+1}$, if $\I(G_{i-1,i})=\varnothing$,
and $q_i\in M_{i,i+1}$ otherwise; 

$3^\circ$~$q_{i+1}q_{i+2}\in M_{i,i+1}$, if $\I(G_{i+1,i+2})=\varnothing$, 
and $q_{i+1}\in M_{i,i+1}$ otherwise. 

If at least one of the parts $G_{i-1,i}$ and $G_{i+1,i+2}$ is empty, we call $M_{i,i+1}$
a {\it boundary cut} of the part~$G_{i,i+1}$.
\end{defin}

\begin{rem}
1) The fact that $M_{i,i+1}\in\mathfrak M^+(G)$ obviously follows from corollary~\ref{l32c1}.

2) Note, that a boundary cut can be not maximal. If~$x$ is the only vertex of the part~$G_{i,i+1}$ adjacent
to~$p$ than the set~$M_{i,i+1}$ can be complemented by an edge~$px$.

3) Also note, that if $M_{i,i+1}\in\mathfrak M(G)$, then $G_{i,i+1}\in\P(M_{i,i+1})$.
\end{rem}

\begin{lem}
\label{lor0} 
Let $F=(p;q_1,\ldots,q_m)$ be a nondegenerate flower. Then the  following statements  hold.

$1)$ If a cutset $Q_{i,i+1}$ can be complemented by an edge~$px$, then $x\in\I(G_{i,i+1})$. 

$2)$ If $M_{i,i+1}\in\mathfrak M(G)$, and the cut $M_{i,i+1}$ can be complemented by an edge, then this edge is~$px$ where $x\in \I(G_{i,i+1})$ is the only vertex of the part~$G_{i,i+1}$  adjacent to~$p$.

\begin{proof}
1) Let $x\not\in\I(G_{i,i+1})$. Then $x\in V(G)\setminus G_{i,i+1}$. Note, that all vertices of the set~$V(G)\setminus G_{i,i+1}$ are in the same connected component of the graph~$G-Q_{i,i+1}$. By lemma~\ref{ll31} we have, that $x$ is  the only vertex of this component which is adjacent to~$p$. From the  other side
the vertex~$p$ must be adjacent to at least one inner vertex of each nonempty part of~$\P(F)$ and, by remark~\ref{r30}, to a common petal of each two neighboring empty parts.

Then there is not more than one nonempty part among all different from~$G_{i,i+1}$ parts. 
If all these parts are empty, then there are at least tree consecutive empty parts, i.e. there at least two
petals not from~$G_{i,i+1}$ adjacent to~$p$, that is impossible.
Thus, there are exactly two nonempty parts in  $\P(F)$ and no neighboring empty parts. That means $m=4$, 
$\I(G_{i-1,i})=\I(G_{i+1,i+2})=\varnothing$ and $\I(G_{i+2,i-1})\ne\varnothing$.
Then $M_{i,i+1}=\{q_iq_{i-1},q_{i+1}q_{i+2},p\}$, and, complementing it by an edge~$px$, we obtain, that~$F$ 
is contained in the cut $M=\{q_iq_{i-1},px,q_{i+1}q_{i+2}\}\in\mathfrak M_3(G)$, i.e. $F$ is a degenerate
flower. We obtain a contradiction.

2) Let the cut $M_{i,i+1}$ can be complemented by an edge~$e$. Then the cutset~$Q_{i,i+1}$ also can be complemented by this edge. Hence, by previous item, $e$ cannot be an edge~$px$, where $x\in\I(G_{i+1,i})$.

Suppose, that $e=q_iv$ (case $e=q_{i+1}v$ is similar). Since $M_{i,i+1}\in\mathfrak M(G)$, we obtain, that $M_{i,i+1}=\{q_i,p,q_{i+1}q_{i+2}\}$. Consider two cases.

\q1. $v\in G_{i+1,i}$. 

Then $v\in G_{i-1,i}$ and $\{q_{i+1},p,q_iv\} \in \mathfrak M_1(G)$, hence, by lemma~\ref{l32} we 
have~ $v=q_{i-1}$ and~$\I(G_{i-1,i})=\varnothing$. But then $q_iv=q_iq_{i-1}\in M_{i,i+1}$.

\q2. $v\in G_{i,i+1}$. 

In this case by lemma~\ref{l30} we obtain, that $\{v,p,q_{i+2}\}$ is a cutset containing~$p$ and dependent with~$Q_{i,i+1}$. Then by lemma~\ref{l300} the flower~$F$ is not maximal, that contradicts  the condition
of lemma.

Thus, the only remained case is $e=px$ where $x\in \I(G_{i,i+1})$. Then by lemma~\ref{ll31}  we have, that~$x$ is the only vertex of the part~$G_{i,i+1}$ adjacent to~$p$.
\end{proof}
\end{lem}

\begin{defin}
If a cut~$M_{i,i+1}$ can be complemented by an edge~$px$ where $x\in \I(G_{i,i+1})$, then denote by~$M^*_{i,i+1}$
the cut, obtained after complementing.
\end{defin}

\begin{figure}
\centerline{\includegraphics{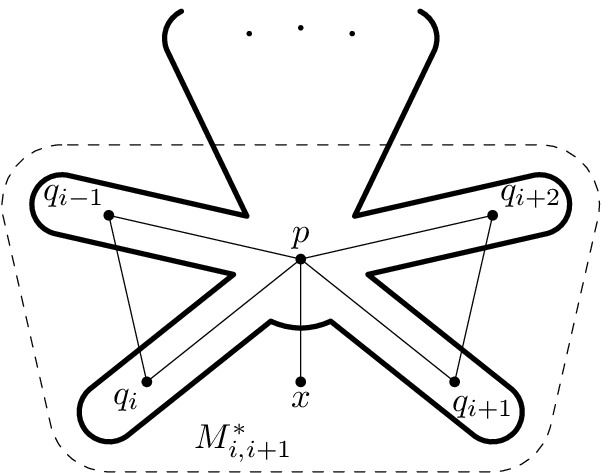}}
\caption{A cut $M^*_{i,i+1}$}\label{ris.5}
\end{figure}

An example of a cut~$M^*_{i,i+1}$ for the case when both parts~$G_{i-1,i}$ and $G_{i+1,i+2}$ are empty is shown on figure~\ref{ris.5}.

\subsection{Sets, splitting basic structures}

In this section we  consider 3-cutsets which split the  vertex set of a cut, a triple cut, or a flower.

\begin{lem}
\label{ll34}
Let maximal nontrivial cut~$M$ and $S\in \mathfrak R_3(G)\setminus\mathfrak R(M)$ be such that
 $S$ splits $V(M)$.
Then $|\P(S)|=2$ and one of the following statements holds.

$1^\circ$ The cut~$M$ is contained in a flower, generated by the set~$\mathfrak S=\mathfrak R(M) \cup \{S,T^M_1,T^M_2\}$.

$2^\circ$ The set~$S$ is contained in a neighborhood of a part of~$\P(M)$  (let it be $G^M_1$). There
exists such an edge~$x_1x_2\in M$ that $S$  separates the vertex~$x_1 \in T^M_1$ from other
vertices of the set~$V(M)\setminus S$, and  $S\setminus G^M_1=\{x_2\}$.

\begin{proof} 
Consider two cases.

\q1. $S\cap \I(G^M_1) \ne\varnothing$ and $S\cap \I(G^M_2) \ne\varnothing$.

In this case the set~$S$ is dependent with all cutsets of~$\mathfrak R(M)$ and also with~$T^M_1$ and~$T^M_2$. 
Thus, the dependence graph of the set~$\mathfrak S$ is connected. 

Prove, that $|S\cap G^M_1|=2$. Let $S\cap G^M_1=\{x\}$. 
Then~$x\in\I(G^M_1)$ and $S\cap T^M_1=\varnothing$.
Whence by corollary~\ref{l1c2} the set~$S$ must separate a single vertex of the set~$T^M_1$
(denote it by~$y$) from other vertices of this set. But then by corollary~\ref{l30c1}
the cut~$M$ can be complemented by an edge~$xy$. We obtain a contradiction with  maximality of the cut~$M$.

Thus, $|S\cap G^M_1|=2$. Similarly, $|S\cap G^M_2|=2$. Therefore, 
$S\cap G^M_1\cap G^M_2\ne\varnothing$. That is there exists a vertex~$p\in M\cap S$.
But this vertex belongs to all cutsets of the set~$\mathfrak S$, whence by corollary~\ref{tr5c1} it follows, that the set~$\mathfrak S$ generates a flower~$F$, which contains the cut~$M$. Since $S\in\mathfrak R(F)$,
then~$|\P(S)|=2$.

\q2. $S\cap \I(G^M_2)=\varnothing$.

In this case~$S\subset\O(G^M_1)$, the set~$S$ is independent with~$T^M_2$ and dependent with~$T^M_1$. 
Similarly to previous item, $|S\cap G^M_1|=2$. Let $S\setminus G^M_1=\{x_2\}$. Since $x_2\in\O(G^M_1)\setminus G^M_1$, then there exists such an edge~$x_1x_2\in M$ that $x_1\in T^M_1$.
Then, clearly, $S$ does not split~$G^M_2$, thus, all vertices of the set~$V(M)\setminus S$, except~$x_1$, 
are in the same connected component of the graph~$G-S$.
Hence~$S$ splits~$T^M_1$ into exactly 2 parts, i.e. $|\P(S)|=2$.
\end{proof}
\end{lem}

\begin{rem}
1) Let us consider in details the second case of the proof of previous lemma
(when $S\setminus G^M_1=\{x_2\}$ and the cutset~$S$ separates~$x_1$ from other vertices of the set~$V(M)$). Obviously, the cutset~$S$ can be complemented by an edge~$x_1x_2$. 
Further two cases are possible: the set~$S\cap T^M_1$ can be either empty or nonempty.

In the first case ($S\cap T^M_1=\varnothing$), clearly,~$x_1x_2$ is a singular  edge.

In the second case let $S\cap T^M_1=\{p\}$. Then the cutsets~$S$, $T^M_1$ and all cutsets of~$\mathfrak R(M)$ containing~$p$ generate a flower. If~$p\in M$, this flower contains~$M$, else  $p$ is an end of an edge~$e\in M$
and our flower contains all vertices of~$V(M)$, except the other end of the edge~$e$.

2) If the cut~$M$ is trivial ($G^M_2=\{x\}$) and $S$ splits~$V(M)$, then~$x\in S$ and either~$S$ separates
one  vertex of the set~$T^M_1$ from two other vertices lying in the same connected component of the 
graph~${G-S}$, or the cut~$M$ is contained in a triple cut with line~$S$.
\end{rem}

Next lemma is about the neighborhood of a triple cut. Remind, that by definition the neighborhood of a triple 
cut~$F=M_1\cup M_2\cup M_3$ is the set~$\O(F)=V(M_1')\cup V(M_2')\cup V(M_3')$, where $M_i'$ is a cut 
from $\mathfrak M_3(G)$ containing~$M_i$ if such  cut exists, and~$M_i'$ coincides with~$M_i$, otherwise.

\begin{lem}
\label{ll35}
Let a triple cut~$F=M_1\cup M_2\cup M_3$ with a line~$S$ and a cutset~$T\in \mathfrak R_3(G)$ be such, that~$T\not\subset\O(F)$ and $T$ splits~$\O(F)$. Then~$|\P(T)|=2$ and the cutset~$T$ is contained in some part~$A_i\in\P(S)$. Moreover, the cutset~$T$ separates a vertex~$x_i \in \I(A_i)$ from other vertices of the 
set~$\O(F)$ and  $T\setminus B_i'=\{x\}$ where $x\in S$, $xx_i\in M_i'$ and $B_i'\in\P(M_i')$ is a part
contained in~$A_i$.

\begin{proof} 
Let~$T$ is dependent with~$S$. Then by lemma~\ref{l3v0} the set~$T$ is subordinated to~$S$, thus,
$T\subset V(F)$. We have a contradiction.

Hence, $T$ is independent with~$S$ and, consequently, $T$ is contained in a part of~$\P(S)$ (let it be~$A_i$). Then, since $T$ splits $\O(F)$ and $T\ne S$, the set~$T$  splits~$V(M_i')$. 
Moreover, $T\not\subset V(M_i')$. If the cut~$M_i'$ is maximal, we apply lemma~\ref{ll34} to $M_i'$ and to the cutset~$T$.  Since the cutsets~$S$ and~$T$ are independent, the statement~$1^\circ$ of lemma~\ref{ll34} 
cannot hold. Hence, the statement~$2^\circ$ of lemma~\ref{ll34} holds, that implies what is to be proved. If the cut~$M_i'$ is not maximal, then, by the definition, $M_i=M_i'\in\mathfrak M_1(G)$ and, since~$T$ is independent with~$S$, the statement of our lemma in this case is clear.
\end{proof}
\end{lem}

\begin{lem}
\label{ll36}
Let a maximal flower $F=(p;q_1,\ldots,q_m)$ and a cutset~$T\in \mathfrak R_3(G)\setminus\mathfrak R(F)$ be such, that $T$ splits $V(F)$. Then~$|\P(T)|=2$ and one of two following statements hold.

$1^\circ$ The set~$T$ separates one vertex of the set~$V(F)$ from other vertices of this set.

$2^\circ$ The set~$T$ separates two neighboring petals~$q_{i+1}$ and~$q_{i+2}$ from other vertices of the 
set~$V(F)$. Moreover,~$\I(G_{i,i+1})=\I(G_{i+2,i+3})=\varnothing$ and $T=\{q_i,x,q_{i+3}\}$, 
where~$x\in\I(G_{i+1,i+2})$ is the only vertex of the part~$G_{i+1,i+2}$ adjacent to~$p$.

\begin{proof} 
Note, that by lemma~\ref{l300} we have~$p\not\in T$. If the cutset~$T$ does not split~$L=\{q_1,\ldots,q_m\}$, 
then~$T$ separates the center~$p$ from all petals of the flower, thus, statement~$1^\circ$ holds. 

Let~$T$ split~$L$. We shall prove that $T\cap L\ne\varnothing$. Indeed, if $\I(G_{i,i+1})\ne \varnothing$ and
 $|T\cap \I(G_{i,i+1})|\le 1$ then by lemma~\ref{l300} the petals~$q_i$ and~$q_{i+1}$ are connected in~${G-T}$. If~$\I(G_{i,i+1})= \varnothing$, then the petals $q_i$ and $q_{i+1}$ are adjacent.
Since there is not more than one such part~$G_{j,j+1}$ that~$|T\cap \I(G_{j,j+1})|\ge 2$, all pairs of neighboring petals (except, maybe, one pair) are not splitted by~ $T$. Hence, if~$T\cap L=\varnothing$, then all
vertices of the set~$L$ are connected in~${G-T}$, we obtain a contradiction.

Note, that~$|T\cap L|\le2$ by lemma~\ref{lor12}. Consider the following two cases.

\q1. $T\cap L=\{q_i\}$.

In this case two other vertices of the cutset~$T$ must be in the same part of~$\P(F)$, otherwise, similarly to proved above, the cutset~$T$ does not split~$L$. Let this part be~$G_{j,j+1}$. Then it is clear, that the cutset~$T$ splits~$L$ into two sets~$\{q_{i+1},\ldots,q_j\}$ and~$\{q_{j+1},\ldots,q_{i-1}\}$. 
Without loss of generality we may assume, that~$p$ and $\{q_{j+1},\ldots,q_{i-1}\}$ lie in the same 
connected component of the graph~${G-T}$. Let us prove, that~$i+1=j$. Indeed, otherwise~$T\cap \I(G_{i,i+2})=\varnothing$, consequently, the cutset~$T$ does not split the part~$G_{i,i+2}$, i.e.~$p$
and~$q_{i+1}$ are connected in~${G-T}$, that contradicts our assumption. Hence, $T$ separates the 
petal~$q_{i+1}=q_j$ from other vertices of~$V(F)$ and statement~$1^\circ$ holds.

\q2. $T\cap L=\{q_i,q_j\}$.

In this case, obviously, the cutset~$T$ splits~$L$ into sets~$\{q_{i+1},\ldots,q_{j-1}\}$ 
and~$\{q_{j+1},\ldots,q_{i-1}\}$. 
Without loss of generality we may assume, that~$p$ and $\{q_{j+1},\ldots,q_{i-1}\}$ lie in the same 
connected component of the graph~${G-T}$. Let us prove, that the other set~$\{q_{i+1},\ldots,q_{j-1}\}$ consists of not more than two vertices. Indeed, otherwise $\I(G_{i,i+2})\cap \I(G_{j-2,j})=\varnothing$ and the 
cutset~$T$ cannot intersect interiors of both these parts. Then, similarly to proved above,
$p$ and $\{q_{i+1},\ldots,q_{j-1}\}$ are connected in~${G-T}$, that contradicts our assumption.

Further, if~$j=i+3$, i.e.~$T$ separates the petals~$q_{i+1},q_{i+2}$ from other vertices of~$V(F)$, then among the parts~$G_{i,i+1}$, $G_{i+1,i+2}$, $G_{i+2,i+3}$ is not more than one nonempty part (because
each nonempty part contains a path, connecting its petal with the center of the flower and consisting of inner vertices of this part). Moreover, among these three parts there are no two neighboring empty parts, because by 
remark~\ref{r30} their common petal is adjacent to the center. It is possible in the only case
$\I(G_{i,i+1})=\I(G_{i+2,i+3})=\varnothing$ and $\I(G_{i+1,i+2})\ne\varnothing$. 
Then, clearly, $T=\{q_i,x,q_{i+3}\}$ where $x\in\I(G_{i+1,i+2})$. In addition, the cutsets~$T$ and $Q_{i+1,i+2}$
are dependent and~$\{p,x\}\in\P(\{T,Q_{i+1,i+2}\})$, hence, $x$ is the only vertex of the part~$G_{i+1,i+2}$
adjacent to~$p$.  Thus,  statement~$2^\circ$ holds.

In all cases it is clear, that~$|\P(T)|=2$, since the cutset~$T$ is dependent with  some cutsets of~$F$ and 
splits each of them into exactly two parts.
\end{proof}
\end{lem}

\begin{rem}
In case~$2^\circ$ of previous lemma it is clear, that the cutset~$T$ is contained in a cut
$M^*_{i,i+1}=\{q_iq_{i+1},px,q_{i+3}q_{i+2}\}$.
\end{rem}

Next lemmas will consider the case~$1^\circ$ of lemma~\ref{ll36} in details.

\begin{lem}
\label{ll37}
Let a maximal flower~$F=(p;q_1,\ldots,q_m)$  and a cutset~$T\in \mathfrak R_3(G)\setminus\mathfrak R(F)$ be such that~$T$ separates a petal~$q_i$ from other vertices of the set~$V(F)$. Then exactly one of two parts
from~$\P(F)$ containing~$q_i$ is empty and the cutset~$T$ consists of the second petal of this part and two vertices of the other part containing~$q_i$.

\begin{proof} 
Let $q_i\in H\in\P(T)$. By the condition,~$q_i$ is the only vertex of~$V(F)$ in  $\I(H)$. 
Then~$\I(H)\cap Q_{i-1,i+1}=\varnothing$, thus, $Q_{i-1,i+1}$ does not split~$H$. 
That means~$H\subset G_{i-1,i+1}$ and, in particular, $T\subset G_{i-1,i+1}$.

Note also, that $p\not\in T$ by lemma~\ref{l300} and, obviously, $q_i\not\in T$. Thus one of the parts~$G_{i-1,i}$ and $G_{i,i+1}$ contains not more than one vertex of the cutset~$T$.
Without loss of generality we may assume, that it is the part~$G_{i-1,i}$. Clearly,~$|T\cap G_{i-1,i}|=1$, otherwise the vertices~$q_i$ and~$q_{i-1}$ are connected in~${G-T}$.
Moreover, if $T\cap G_{i-1,i}\ne\{q_i\}$, then by lemma~\ref{l300} there is a path connecting~$q_i$ 
and~$q_{i-1}$ in the part~$G_{i-1,i}$ which do not intersect~$T$. This is impossible.
Thus, $T\cap G_{i-1,i}=\{q_{i-1}\}$, i.e. $T\cap\I(G_{i-1,i})=\varnothing$. 
But~$T$  separates from each other the vertices~$p,q_i\in G_{i-1,i}$.  It is possible only if~$\I(G_{i-1,i})=\varnothing$. 
Then~$\I(G_{i,i+1})\ne\varnothing$, since otherwise by remark~\ref{r30} the vertices $q_i$ and~$p$ are adjacent.
\end{proof}
\end{lem}

\begin{lem}
\label{ll38}
Let a cutset~$T\in \mathfrak R_3(G)$ separate a center of nondegenerate flower~$F=(p;q_1,\ldots,q_m)$ from
other vertices of the set~$V(F)$.  Then~$T$ separates~$p$ from
other vertices of the graph~$G$.  Moreover, $m\le6$, there are not more than three nonempty parts in~$\P(F)$, the interior of every nonempty part contains exactly one vertex of the set~$T$, 
and the boundary of every nonempty part does not intersect~$T$.

\begin{proof} 
Let $\P(F)$ contain  $k$ nonempty parts and  $\ell$ empty parts.
If~$\I(G_{i,i+1})\ne\varnothing$, then~$Q_{i,i+1}$ is a cutset dependent with~$T$.
Consequently, $T\cap \I(G_{i,i+1})\ne\varnothing$. Thus, $k\le3$.
Further, if $\I(G_{j-1,j})=\I(G_{j,j+1})=\varnothing$, then by remark~\ref{r30} the vertices~$p$ and~$q_j$ 
are adjacent, hence, $q_j\in T$. Note, that empty parts of~$\P(F)$ are divided into not more than~$k$ sequences,
which give us at least~$\ell-k$ petals adjacent with~$p$. Thus, $\ell=k+(\ell-k)\le3$, hence, $m=k+\ell\le6$.

Let $|T\cap G_{i,i+1}|=2$.
Then~$|T\cap \I(G_{i+1,i})|=1$ and by lemma~\ref{ll31} the cutset $Q_{i,i+1}$ can be complemented by an edge~$px$, where
$T\cap \I(G_{i+1,i})=\{x\}$. This contradicts lemma~\ref{lor0}.

Thus, $|T\cap G_{i,i+1}|\le1$ for every part~$G_{i,i+1}$. Hence, if $\I(G_{i,i+1})\ne\varnothing$, 
then $|T\cap \I(G_{i,i+1})|=1$ and~$T\cap Q_{i,i+1}=\varnothing$. 
In addition, if~$T\cap \I(G_{i,i+1})=\{u\}$, then~ $\{p,u\}\in\P(\{T,Q_{i,i+1}\})$ by corollary~\ref{l1c2}, i.e. the cutset~$T$ separates~$p$ from other vertices of the part~$G_{i,i+1}$. Since this condition holds for every nonempty part, then~$T$ separates~$p$ from other vertices of the  graph~$G$.
\end{proof}
\end{lem}

\subsection{Singular flowers. The neighborhood of a flower}

\begin{defin}
\label{dnc1}
We call a flower~$F=(p;q_1,\ldots,q_m)$ {\it singular}, if $d(p)=3$ and {\it nonsingular} otherwise.
Let  {\it neighborhood} of the center of this singular flower be the set~$T(p)$ consisting of all adjacent to~$p$ vertices. 
\end{defin}

Note, that if $p$ is the center of a singular flower, then $T(p)\in\mathfrak R_3(G)$.
Moreover, by lemma~\ref{ll38} the interior of each nonempty part of~$\P(F)$ contains exactly one vertex of~$T(p)$, and its boundary does not intersect~$T(p)$.

\begin{defin}
\label{dor}
Let $F=(p;q_1,\ldots,q_m)$ be a maximal nondegenerate flower and~$G_{i,i+1}\in\P(F)$. 

If the flower~$F$ is singular and the part~$G_{i,i+1}$ is nonempty, then denote by~$u_{i,i+1}$ 
the only vertex of~$G_{i,i+1}\cap T(p)$.

If the flower~$F$ is nonsingular, there is exactly one vertex adjacent to~$p$ in the part~$G_{i,i+1}$
and $\I(G_{i-1,i})=\I(G_{i+1,i+2})=\varnothing$. Then also denote by $u_{i,i+1}$
the only adjacent to~$p$ vertex of~$G_{i,i+1}$. 

In all other cases we set~$u_{i,i+1}=p$.

The set $\O(F)=V(F)\cup \{u_{1,2},u_{2,3},\ldots,u_{m,1}\}$ we call the {\it neighborhood} of the flower~$F$.
\end{defin}

\begin{rem}
Note, that if $u_{i,i+1}\ne p$, then the set~$M_{i,i+1}$ can be complemented by an edge~$pu_{i,i+1}$, i.e. $pu_{i,i+1}\in M^*_{i,i+1}$. 

If $F$ is a maximal nondegenerate singular flower, then by definition~$\O(F)=V(F)\cup T(p)$.

If~$F$ is a maximal nondegenerate  nonsingular flower and $u_{i,i+1}\ne p$, then by definition
$\I(G_{i-1,i})=\I(G_{i+1,i+2})=\varnothing$, consequently, 
$M^*_{i,i+1}=\{q_{i-1}q_i,pu_{i,i+1},q_{i+2}q_{i+1}\}\in\mathfrak M_3(G)$.
On the other side, if~$M^*_{i,i+1}\in\mathfrak M_3(G)$, then, clearly,~$u_{i,i+1}\ne p$.
\end{rem}

\begin{defin}
Let $u_{i,i+1}\ne p$. Set the notations $M'_{i,i+1}=M^*_{i,i+1}$ and 
$Q'_{i,i+1}=\{q_i,u_{i,i+1},q_{i+1}\}$. If~ $u_{i,i+1}=p$, we set~$M'_{i,i+1}=M_{i,i+1}$ and~$Q'_{i,i+1}=Q_{i,i+1}$.  

Let  $G'_{i,i+1}=G_{i,i+1}\setminus \{p\}$ if either $u_{i,i+1}\ne p$, or
$\I(G_{i,i+1})=\varnothing$ and at least one of the vertices $u_{i-1,i}$ and $u_{i+1,i+2}$ differs from~$p$. 
In all other cases we set $G'_{i,i+1}=G_{i,i+1}$.

If $M'_{i,i+1}\in\mathfrak M(G)$, then denote by $\O(G'_{i,i+1})$ the neighborhood of~$G'_{i,i+1}$ as of a part of~$\P(M'_{i,i+1})$.
We call the cut~$M'_{i,i+1}$  {\it boundary} cut of the part~$G'_{i,i+1}$.
If $M'_{i,i+1}\in\mathfrak M_0(G)$ we set $\O(G'_{i,i+1})=G'_{i,i+1}$.
\end{defin}

For $u_{i,i+1}\ne p$ it is easy to see that~$G'_{i,i+1}$ is a part of~$\P(M'_{i,i+1})$ contained
in~$G_{i,i+1}$ and $\I(G'_{i,i+1})=\I(G_{i,i+1})\setminus \{u_{i,i+1}\}$.

\begin{lem}
\label{ll39}
Let $F=(p;q_1,\ldots,q_m)$ be a maximal nondegenerate flower. Then the following statements hold.

$1)$ If $T\in\mathfrak R_3(G)$ and $T\subset\O(F)$, then either $T\in \mathfrak R(F)$,  or~$T$ is contained in~$M_{i,i+1}'$ for some~$i$,  or~$T=T(p)$  (the last is possible only for a singular flower~$F$). 
All sets described above, except sets~$Q'_{i,i+1}$, are  cutsets splitting the graph~$G$
into exactly two parts. Each of these cutsets splits~$\O(F)$.
The set~$Q'_{i,i+1}$ does not split~$\O(F)$  and is a cutset if and only if~$\I(G'_{i,i+1})\ne\varnothing$.

$2)$ If a cutset~$S\in\mathfrak R_3(G)$ splits~$\O(F)$ and $S\not\subset\O(F)$, then
$|\P(S)|=2$ and $S$ separates one vertex of~$\O(F)$ from the other vertices of this set. The vertex separated by~$S$ is not the center of~$F$.

\begin{proof} 
1) If~$T\subset V(F)$, then by corollary~\ref{lor12c1} we have, that $T$ is a set of the flower~$F$. 
Then either~$T\in\mathfrak R(F)$, or~$T$ is a boundary of a nonempty part~$G_{i,i+1}\in\P(F)$ 
and is contained in~$M_{i,i+1}'$.

Let $T\not\subset V(F)$. Then by definition of the neighborhood of a flower there exists such~$i$,  that
 $u_{i,i+1}\in T$ and $u_{i,i+1}\ne p$. Note, that if~$T$ does not split~$V(F)$, then
$T\subset G_{i,i+1}\cap\O(F)=\{p,q_i,q_{i+1},u_{i,i+1}\}\subset V(M'_{i,i+1})$.
By corollary~\ref{l31c1} we obtain, that~$T$ is contained in~$M'_{i,i+1}$.

The only remaining case is~$T\not\subset V(F)$ and~$T$ splits $V(F)$.
By lemma~\ref{ll36} there are 3 possible subcases.

\q1. {\it $T$ separates $p$ from other vertices of~$V(F)$.} Then by lemma~\ref{ll38} the flower~$F$ 
is singular and~$T=T(p)$.

\q2. {\it $T$ separates one petal  of the flower~$F$ from other vertices of~$V(F)$.}  
Since $T\cap I(G_{i,i+1})\ne\varnothing$, then this petal is~$q_i$ or~$q_{i+1}$.  
Then by lemma~\ref{ll37} two vertices of~$T$ belong to
$G_{i,i+1}\cap\O(F)\subset V(M'_{i,i+1})$, and the third vertex belongs to a neighboring with~$G_{i,i+1}$
empty part of~$\P(F)$, i.e. also belongs to~$V(M'_{i,i+1})$. Thus, 
$T\subset V(M'_{i,i+1})$, hence, $T$ is contained in $M'_{i,i+1}$.

\q3. {\it $T$ separates petals~$q_i$ and~$q_{i+1}$ from other vertices of~$V(F)$.}
Thus, by lemma~\ref{ll36} we have $\I(G_{i-1,i})=\I(G_{i+1,i+2})=\varnothing$ and
$T=\{q_{i-1},u_{i,i+1},q_{i+2}\}\subset M'_{i,i+1}$.

\smallskip
By the properties of a flower, every its inner set is a cutset, which splits the graph~$G$  into two parts
and splits~$V(F)$ (consequently, it also splits~$\O(F)$).
By lemma~\ref{l30} the same statement holds for every inner set of the cut~$M'_{i,i+1}$,
and by lemma~\ref{ll38}~--- for the set~$T(P)$ (if the flower~$F$ is singular).
If $Q_{i,i+1}\ne Q'_{i,i+1}$, then, obviously, $Q_{i,i+1}$ is a cutset separating~ $u_{i,i+1}$ from other
vertices of the set~$\O(F)$ and, consequently, splits the graph~$G$  into exactly two parts. All remaining sets are of type~$Q'_{i,i+1}$. Clearly,~$Q'_{i,i+1}$ does not split~$\O(F)$ and is a cutset if and only if~$\I(G'_{i,i+1})\ne\varnothing$.

\smallskip
2) If $S$ splits $V(F)$, then by lemma~\ref{l300} we have~$p\not\in S$.
Further, by lemma~\ref{ll36}, the cutset~$S$ splits the graph into exactly two parts and
separates not more than two vertices of the set~$V(F)$ from other vertices of this set. If~$S$ separates
two vertices, then~$S\subset\O(F)$, that contradicts the condition. In addition, by lemma~\ref{ll38} the only 
3-cutset separating~$p$ from other vertices of the set~$V(F)$ is~$T(p)$. But $T(p)\subset \O(F)$.

Hence, $S$ separates one petal of the flower from  other vertices of the set~$V(F)$.
In addition, since~$p\not\in S$, then all vertices~$u_{i,i+1}$ and~$p$ belong to the same connected component of the graph~$G-S$. Thus, $S$ separates one petal of the flower from other vertices of the set~$\O(F)$.

If~$S$ does not split~$V(F)$, then $S$ is independent with all sets of the flower~$F$, i.e.~$S$ is contained in some part~$G_{i,i+1}$. But then~$S$ can separate from the other vertices of the set~$\O(F)$  not more than 
one vertex~$u_{i,i+1}$, and it is possible only in the case~$u_{i,i+1}\ne p$. Moreover, by
lemma~\ref{lds1} the cutset $S$ does not split~$G_{i+1,i}$.
In this case, clearly,  $p\in S$, and every part of~$\P(S)$, which is a subset of~$G_{i,i+1}$ contains a vertex
adjacent to~$p$. But, since $u_{i,i+1}\ne p$, there is only one such vertex in~$G_{i,i+1}$, 
consequently,  only one part of $\P(S)$ is contained in~$G_{i,i+1}$. On the other side, by lemma~\ref{lds1}, 
the cutset $S$ does not split~$G_{i+1,i}$, hence, there is only one part of~$\P(S)$ not contained in~$G_{i,i+1}$ i.e.~$|\P(S)|=2$. 
\end{proof}
\end{lem}

\begin{rem}
\label{r31}
Let us describe nondegenerate singular flower $F(p;q_1,\ldots,q_m)$ in details.
Consider  several cases.

1) Let~$\P(F)$ contain three nonempty parts. Clearly,~$T(p)$ contains a vertex from the interior of each nonempty part, hence,~$p$ is not adjacent to any petal of~$F$.
Whence by remark~\ref{r30} it follows, that no two empty parts of~$\P(F)$ are neighboring. 
Thus,~$\P(F)$ contains 4, 5 or 6 parts. A singular flower with 6 parts and three nonempty parts among them
is shown on figure~\ref{ris.6}.

\begin{figure}
\centerline{\includegraphics{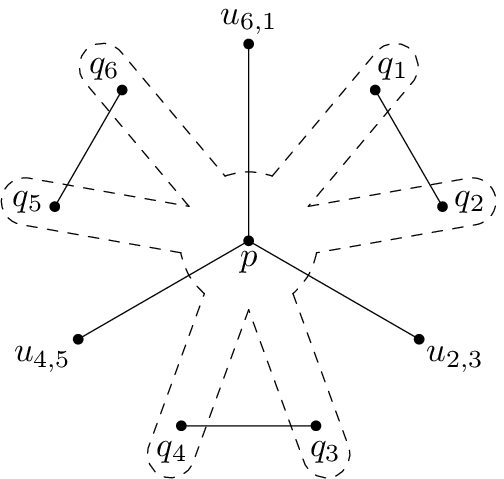}}
\caption{A singular flower with six parts and three nonempty parts}\label{ris.6}
\end{figure}

2) Let~$\P(F)$ contain two nonempty parts, then the interior of each nonempty part contains exactly one vertex of the set~$T(p)$, and the third vertex of~$T(p)$ is a petal~$q_i$, not belonging to any nonempty part of~$\P(F)$ (i.e., the parts~$G_{i-1,i}$ and~$G_{i,i+1}$ are empty). Since all vertices 
adjacent to~$p$ belong to~$T(p)$, the decomposition~$\P(F)$ except four parts mentioned above
can contain not more than one empty part. Moreover, this empty part  must be neighboring with two nonempty parts.
Thus~$\P(F)$ consists of four or five parts.

If $|\P(F)|=4$, then this parts can be enumerated such that $G_{1,2}\ni u$ and $G_{4,1}\ni v$ are nonempty parts, $G_{2,3}$ and~$G_{3,4}$ are empty parts, and $T(p)=\{u,v,q_3\}$.
Note, that~$F'=(q_3;q_1,u,p,v)$ is a flower of the same type, its center~$q_3$ is separated from other vertices of the graph~$G$ by the set  $\{q_2,p,q_4\}$. It is easy to see, that~$\O(F)=\O(F')$.

\begin{figure}
\centerline{\includegraphics{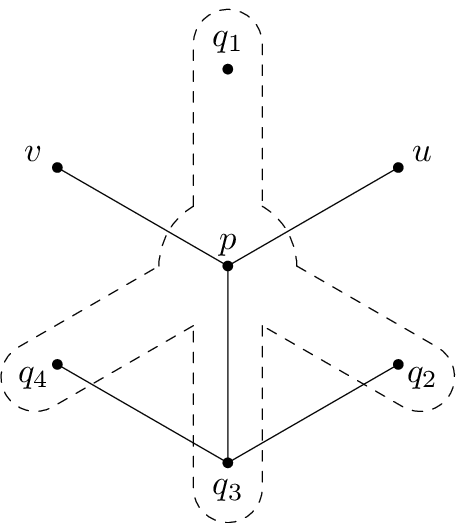}\qquad\qquad\qquad\includegraphics{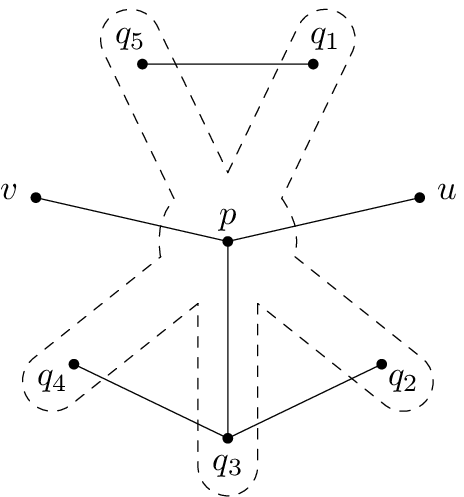}}
\caption{Singular petals with 4 or 5 parts and 2 nonempty parts}\label{ris.7}
\end{figure}

If $|\P(F)|=4$, then the parts can be enumerated such, that  $G_{1,2}\ni u$ and 
$G_{4,5}\ni v$ are nonempty parts, $G_{2,3}$, $G_{3,4}$ and~$G_{5,1}$ are empty parts, and $T(p)=\{u,v,q_3\}$.

Note, that in this case~$F'=(q_3;q_1,u,p,v,q_5)$ is a flower of the same type, its center~$q_3$ is separated from other vertices of the graph~$G$ by the set  $\{q_2,p,q_4\}$. It is easy to see, that~$\O(F)=\O(F')$. 

Singular flowers with 4 or 5 parts and 2 nonempty parts are shown on the figure~\ref{ris.7}.

3) The case when $\P(F)$ contains exactly one nonempty part is impossible. Indeed, if it happens, the interior of this nonempty part contains a vertex~$u\in T(p)$, and two other vertices of~$T(p)$ are petals of our flower.
There are at least three other parts of~$\P(F)$, all of them  are empty. By remark~\ref{r30} the common petal of two empty  parts is adjacent to the center~$p$, thus, this petal belongs to~$T(p)$. Hence, there are exactly three empty parts in~$\P(F)$ --- let them be~$G_{1,2}$, $G_{2,3}$ and $G_{3,4}$. 
Now it is easy to see, that $\{up,q_1q_2,q_4q_3\}\in \mathfrak M_3(G)$ is a cut containing~$F$. That
contradicts maximality of~$F$.
\end{rem}

\begin{lem}
\label{ll391}
Let~$F$ be a singular nondegenerate flower. Then the following statements hold.

$1)$ For every part~$G_{i,i+1}\in\P(F)$ we have $G_{i,i+1}'\ne G_{i,i+1}$.

$2)$ The set~$T(p)$ consists of all vertices~$u_{i,i+1}$ different from~$p$, 
and all vertices~$q_j$ for which  $\I(G_{j-1,j})=\I(G_{j,j+1})=\varnothing$, 
$u_{j-2,j-1}\ne p$ and $u_{j+1,j+2}\ne p$.

\begin{proof}

1) Note, that in a singular flower~$u_{i,i+1}\ne p$ (i.e. $G_{i,i+1}'\ne G_{i,i+1}$) if and only if~$\I(G_{i,i+1})\ne\varnothing$. 
Moreover, it follows from the classification of singular flowers (see remark~\ref{r31}), that for each empty part of~$\P(F)$ there is a nonempty neighboring part of~$\P(F)$. Thus, for empty parts of~$\P(F)$ we also have~$G_{i,i+1}'\ne G_{i,i+1}$.

2) Note, that $T(p)$ consists of different from~$p$  vertices~$u_{i,i+1}$ and all petals~$q_i$ belonging to two empty parts. It follows from remark~\ref{r31}, that there are no three empty consecutive parts in~$\P(F)$, thus,  $q_i\in T(p)$ means, that $I(G_{i-1,i})=I(G_{i,i+1})=\varnothing$, $I(G_{i-2,i-1})\ne\varnothing$ and $I(G_{i+1,i+2})\ne\varnothing$. Then by item~1 we have $u_{i-2,i-1}\ne p$ and~$u_{i+1,i+2}\ne p$.
\end{proof}
\end{lem}

\subsection{A connection between triple cuts and other basic structures}

A line of triple cut splits a graph into three parts and, consequently, it cannot be an inner set of a flower
or a cut. Thus, the vertex set of a triple cut cannot be a subset of a vertex set of a flower or of a cut.

Moreover, it follows from lemmas~\ref{ll34} and~\ref{ll39}, that a line of a triple cut cannot split a vertex 
set of a nontrivial cut  or a neighborhood of a flower.
Thus, if the vertex set of a nontrivial cut or of a flower is contained in the vertex set of a triple cut~$F=M_1\cup M_2\cup M_3$ then it is contained in a vertex set of one of  three cuts $M_1$, $M_2$, $M_3$.
If a vertex set of a nontrivial cut or of a flower is contained in $\O(F)$, then  it is contained in the vertex set of one of three cuts $M'_1$, $M'_2$, $M'_3$.
Note also, that edges connecting a vertex of degree 3 belonging to the line of a triple cut with three vertices of its neighborhood form a trivial cut which is contained in our triple cut.

\begin{defin}
We say, that a cut or a flower  {\it is contained} in a triple cut, if its vertex set is contained in a vertex set of this triple cut.

We say, that a cut or a flower  {\it is contained} in a neighborhood of a  triple cut, if its vertex set is contained in this neighborhood.
\end{defin}

\section{Complexes}

We represent the set~$\mathfrak R_3(G)$ as a union of several subsets --- structural units of decomposition,
which are constructed on base of structures described above.
We call these subsets  {\it complexes}. 

In this section we present all types of complexes, describe all 3-cutsets of each complex and the decomposition of the graph~$G$ by cutsets of one complex. Further  with the help of theorem of decomposition~\cite{k05} we construct a hypertree of  relative  position of different complexes. 
As a result we obtain a full description of relative position of all 3-cutsets in a triconnected graph.

\begin{defin}
We call a cutset $S\in\mathfrak R_3(G)$ {\it single}, if it is independent with any other cutset 
of~$\mathfrak R_3(G)$. Otherwise, we call cutset $S$  {\it nonsingle}. 
\end{defin}

A {\it single complex} is a complex consisting of one single cutset. Further we describe other 
types of complexes.

\subsection{Triple complexes}

\begin{defin}
For any triple cut~$F$  let the  set consisting of all 3-cutsets contained in $\O(F)$, except boundaries
of the neighborhood of~$F$  be a {\it triple complex}. Let the line of $F$ be  the {\it line} of this triple complex, and  boundaries of the neighborhood of $F$ be {\it boundaries} of this triple complex.
\end{defin}

Let $F=M_1\cup M_2\cup M_3$ be a triple cut with line~$S$, and $\O(F)=V(M_1')\cup V(M_2')\cup V(M_3')$
be its neighborhood. Let $\P(S)=\{A_1,A_2,A_3\}$, and parts $B_i\in\P(M_i)$ and   $B_i'\in\P(M_i')$ are such that  $B'_i\subset B_i\subset A_i$. (All these parts are discussed in details in the Section~2.3.)

By lemma~\ref{ll11}, the triple complex~$\mathfrak C(F)$ consists of~$S$, all trivial cutsets subordinated to~$S$ (there are not more than three such cutsets) and inner cutsets of the cuts~$M_1'$, $M_2'$, $M_3'$.

Let us describe all parts of~$\P(\mathfrak C(F))$. If all cuts~$M_i'$ are nontrivial, 
then~$\P(\mathfrak C(F))$ consists of small parts~$\{x,x_i\}$ (where $x\in S$ and $xx_i\in M'_i$) and 
normal parts~$B_i'$. 

If the cut~$M_i'$ is trivial, then~$|B_i'|=1$. Let~$B_i'=\{y\}$. Then the part~$B_i$ consists of the vertex~$y$ 
and those vertices of the set~$S$ which degree is more, than three. If there are no such vertices,
then all parts of~$\P(\mathfrak C(F))$, contained in~$A_i$ are parts~$\{x,y\}$, where $x\in S$.

Let~$S$ contains a  vertex of degree more, than 3. Then~$B_i\in\P(\mathfrak C(F))$. Moreover, the 
part~$B_i$ is small, if exactly two vertices of~$S$ has degree 3. Otherwise, there is exactly one vertex of degree three in $S$, in this case $B_i$ is a normal part.

It follows from lemma~\ref{ll35}, that for any set~$R\in\mathfrak R_3(G)\setminus \mathfrak C(F)$ 
there exists a unique nonempty part~$A\in \P(\mathfrak C(F))$ such that $R\subset \O(A)$ 
and either $R=\R(A)$, or $R\cap\I(A)\ne\varnothing$.

\subsection{A complex of nondegenerate flower}

Let there exists  such a flower~$F$ in the graph~$G$ that all  parts of~$\P(F)$ are empty. Then all vertices of the graph~$G$ are vertices of~$F$, each petal is adjacent to the center and   two neighboring petals. 
All  3-cutsets of the graph~$G$ are sets of the flower~$F$. Note, that in this case the graph~$G$ is a ``wheel'' (see~\cite{T2}). 
Further we assume, that for every  flower~$F$ in the graph~$G$ the decomposition~$\P(F)$ contains a nonempty part.

In this section we consider a maximal nondegenerate flower~$F=(p;q_1,q_2,\ldots,q_m)$. As it was shown above, there are two essentially different cases: the flower~$F$ can be singular of nonsingular.

\begin{defin}
Let  a {\it complex} $\mathfrak C(F)$ of a flower~$F$ be the set of all 3-cutsets contained in its neighborhood~$\O(F)$ which split~$\O(F)$. 
We call {\it boundaries} of the complex $\mathfrak C(F)$  boundaries of all normal parts of~$\P(\mathfrak C(F))$.
\end{defin}

It follows from lemma~\ref{ll39}, that $\mathfrak C(F)$ consists of cutsets of~$\mathfrak R(F)$, cutsets contained in 
cuts~$M'_{i,i+1}$ and not coinciding with~$Q_{i,i+1}'$ (here it is enough to consider only such~$i$ for which~$u_{i,i+1}\ne p$) and, if the flower~$F$ is singular, the cutset~ $T(p)$. Boundaries of the complex of nondegenerate flower do not belong to this complex!

\begin{lem}
\label{l40}
Let $G_{i,i+1}\in\P(\mathfrak R(F))$ and $G_{i,i+1}'\ne G_{i,i+1}$. Then there exists a 
cutset~$S_{i,i+1}\in \mathfrak C(F)$ separating $p$ from $G_{i,i+1}'$. Moreover,
$S_{i,i+1}\cap G_{i,i+1}=\{u_{i,i+1}\}$ if $\I(G_{i,i+1})\ne\varnothing$ and
$S_{i,i+1}\cap G_{i,i+1}$ consists of one of the petals~$q_i$ and $q_{i+1}$, otherwise.

\begin{proof}
If $\I(G_{i,i+1})=\varnothing$, then, since $G_{i,i+1}'\ne G_{i,i+1}$, at least
one of the vertices~$u_{i-1,i}$ and~$u_{i+1,i+2}$ does not coincide with~$p$. Without loss of generality, let it be $u_{i-1,i}$. Then the set  $S_{i,i+1}=\{q_{i-1},u_{i-1,i},q_{i+1}\}$ is what we want. If~$\I(G_{i,i+1})\ne\varnothing$, then~$S_{i,i+1}=T(p)$ for a singular flower
and~$S_{i,i+1}=\{q_{i-1},u_{i,i+1},q_{i+2}\}$ for a nonsingular flower  is the desired set.
\end{proof}
\end{lem}

To describe  parts of decomposition of the graph by a  complex of  flower we need to generalize
 the notion of the neighborhood of the center for the case of nonsingular flower.
For this purpose we get help of the property of the neighborhood of the center of a singular flower, proved in
item~2 of lemma~\ref{ll391}.

\begin{defin}
\label{dnc2}
Let the {\it neighborhood} of the center of a flower~$F$ be  a set~$T(p)$ consisting of all vertices
$u_{i,i+1}$ different from~$p$, and all petals~$q_j$ for which
$\I(G_{j-1,j})=\I(G_{j,j+1})=\varnothing$, $u_{j-2,j-1}\ne p$ and~$u_{j+1,j+2}\ne p$.
\end{defin}

\begin{figure}
\centerline{\includegraphics{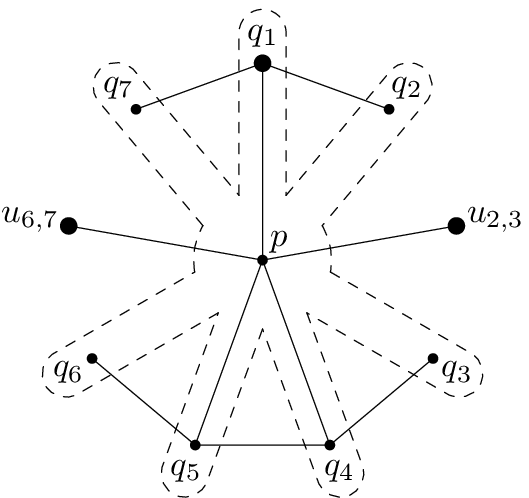}}
\caption{The neighborhood of the center of nonsingular flower}\label{ris.8}
\end{figure}

For example, the neighborhood of the center of the  flower shown on figure~\ref{ris.8} consists of
vertices~$q_1$, $u_{2,3}$ and $u_{6,7}$. The petals~$q_4$ and~$q_5$ do not belong to the neighborhood, 
since~$u_{5,6}=u_{3,4}=p$.

\begin{lem}
\label{l41}
The set $\P(\mathfrak C(F))$ consists of all parts~$G_{i,i+1}'$ and small parts~$\{p,x\}$, 
where~$x\in T(p)$. A part~$G_{i,i+1}'$ is small if and only if~$\I(G_{i,i+1})=\varnothing$ and at least one of the vertices~$u_{i-1,i}$  and~$u_{i+1,i+2}$ is different from~$p$. In addition, if the part~$G_{i,i+1}'$ is normal, then~$\R(G_{i,i+1}')=Q_{i,i+1}'$.

\begin{proof}
By the definition $\P(\mathfrak R(F))=\P(F)=\{G_{1,2},G_{2,3},\ldots,G_{m,1}\}$. 
Let us see, how other sets of~$\mathfrak C(F)$ can split the  parts of~$\P(F)$.

At first we shall prove, that $G_{i,i+1}'\in\P(\mathfrak C(F))$ for all $i$.
Note, that if $G_{i,i+1}'\ne G_{i,i+1}$, then by lemma~\ref{l40} there exists a 
cutset $S_{i,i+1}\in \mathfrak C(F)$ separating~$p$ from~$G_{i,i+1}'$. Thus, it is enough to prove, that
no set of~$\mathfrak C(F)$ splits~$G_{i,i+1}'$. Consider several cases.

\q1. Let $G_{i,i+1}'=G_{i,i+1}$. Then by lemma~\ref{ll391} the flower~$F$ is nonsingular.

\q{1.1}. Let $\I(G_{i,i+1})=\varnothing$. Then $u_{i-1,i}=u_{i+1,i+2}=p$ and~$q_i$ is adjacent to~$q_{i+1}$. Hence,  the cutset~$T$, splitting~$G_{i,i+1}'$, must separate~$p$ from~$\{q_i,q_{i+1}\}$. Without loss of generality we may assume, that $q_i\not\in T$. Thus, by lemmas~\ref{ll36} and~\ref{ll37} we have $T\cap\I(G_{i-1,i})\ne\varnothing$, whence $T\not\in \mathfrak C(F)$, since $u_{i-1,i}=p$.

\q{1.2}. Let $\I(G_{i,i+1})\ne\varnothing$. Since $u_{i,i+1}=p$,  then no cutset of~$\mathfrak C(F)$  intersect~$\I(G_{i,i+1})$ or coincide with~$Q_{i,i+1}$. But $G_{i,i+1}$ is a union of several parts of~$\P(Q_{i,i+1})$. Hence, by lemma~\ref{lds1}, no cutset of~$\mathfrak C(F)$ can split~$G_{i,i+1}=G_{i,i+1}'$.

\q2. Let $G_{i,i+1}'\ne G_{i,i+1}$. The case~$\I(G_{i,i+1})=\varnothing$ is clear, since in this case the vertices~$q_i$ and~$q_{i+1}$ are adjacent and no cutset can split~$G'_{i,i+1}=\{q_i,q_{i+1}\}$.
Note, that it is the  only case when the part~$G_{i,i+1}'$ is small.
Thus, it is enough to consider the case~$\I(G_{i,i+1})\ne\varnothing$. We divide it into two subcases.

\q{2.1}. If $\I(G_{i,i+1}')\ne\varnothing$, then the part~$G_{i,i+1}'$ is a union of several parts of~$\P(Q_{i,i+1}')$.  Since no cutset of~$\mathfrak C(F)$  intersect~$\I(G'_{i,i+1})$ or coincide with~$Q'_{i,i+1}$, then  by lemma~\ref{lds1} no set of~$\mathfrak C(F)$ can split~$G_{i,i+1}'$.

\q{2.2}. If $\I(G_{i,i+1}')=\varnothing$, then $\I(G_{i,i+1})=\{u_{i,i+1}\}$, i.e. the vertex~$u_{i,i+1}$ 
is adjacent to~$q_i$ and~$q_{i+1}$. Further, by lemma~\ref{l300} there exist a path between~$q_i$ and~$q_{i+1}$,
avoiding~$u_{i,i+1}$, which  inner  vertices belong to~$\I(G_{i,i+1})$. Hence, the vertices~$q_i$ and~$q_{i+1}$ are also adjacent and no cutset can split~$G_{i,i+1}'=\{q_i,u_{i,i+1},q_{i+1}\}$.

\smallskip
Now all the sets~$G_{i,i+1}'$ belong to $\P(\mathfrak C(F))$. 
It follows from the definition of~$G_{i,i+1}'$ and~$Q_{i,i+1}'$, that~$\R(G_{i,i+1}')=Q_{i,i+1}'$ if the 
part~$G_{i,i+1}'$ is normal.

Note, that all the sets~$\{p,x\}$ where $x\in T(p)$ also belong to~$\P(\mathfrak C(F))$.
Indeed, $p$ and $x$ are adjacent and $\{p,x\}$ can be separated from other vertices of the graph~$G$ 
by several cutsets of~$\mathfrak C(F)$. If~$x=u_{i,i+1}$, then the sets~$Q_{i,i+1}$ and~$S_{i,i+1}$ (which was constructed in lemma~\ref{l40}) fit for this purpose. Otherwise, if~$x=q_j$, then we use the sets~$Q_{j-1,j+1}$, $S_{i-1,i}$  and~$S_{i,i+1}$. 

Prove, that there are no  other parts. Let $H\in\P(\mathfrak C(F))$ be another part. 
Clearly, the part $H$ is contained in one of the sets~$G_{i,i+1}$ and~$H\not\subset G_{i,i+1}'$, hence, $G_{i,i+1}'\ne G_{i,i+1}$ and~$p\in H$. 

If~$\I(G_{i,i+1})\ne\varnothing$, then $u_{i,i+1}\ne p$ and~$H=\{p,u_{i,i+1}\}$, since by  lemma~\ref{l40} there exists a set~$S_{i,i+1}\in \mathfrak C(F)$ separating~$p$ from other vertices  of~$G_{i,i+1}$. 

If~$\I(G_{i,i+1})=\varnothing$, then either~$H=\{p,q_i\}$, or $H=\{p,q_{i+1}\}$. Without loss of generality we assume, that~$H=\{p,q_i\}$. Then~$H\subset G_{i-1,i}$, hence, $\I(G_{i-1,i})=\varnothing$ (otherwise we get help
of proved above), i.e. the vertices~$p$ and~$q_i$ are adjacent. 
Further, by lemma~\ref{l40} there exist sets~$S_{i-1,i}$ and~$S_{i,i+1}$, separating~$p$ from~$q_{i-1}$ 
and~$q_{i+1}$ respectively. Then by lemmas~\ref{ll36} and~\ref{ll37} we have~$S_{i-1,i}\cap\I(G_{i-2,i-1})\ne\varnothing$ and~$S_{i,i+1}\cap\I(G_{i+1,i+2})\ne\varnothing$, 
hence,~$u_{i-2,i-1}\ne p$ and~$u_{i+1,i+2}\ne p$. That means $q_i\in T(p)$.
\end{proof}
\end{lem}

For every normal part of~$\P(\mathfrak C(F))$ we define its neighborhood.

\begin{defin}
For every normal part~$G_{i,i+1}'\in\P(F)$ let its neighborhood be the set~$\O(G_{i,i+1}')=G_{i,i+1}'\cup V(M_{i,i+1}')$.
\end{defin}

\begin{rem}
If $M_{i,i+1}'=Q_{i,i+1}$, then~$\O(G_{i,i+1}')=G_{i,i+1}'=G_{i,i+1}$.
In all other cases the normal part~$G_{i,i+1}'\in\P(\mathfrak C(F))$ has neighboring cut~$M_{i,i+1}'$. 
Then the neighborhoods of~$G_{i,i+1}'$ as a part of~$\P(\mathfrak C(F))$ and as a part of~$\P(M_{i,i+1}')$ coincide.
\end{rem}

\begin{thm}
\label{tcr} 
Let $F=(p;q_1,\ldots,q_m)$ be a maximal nondegenerate flower,~$R\in\mathfrak R_3(G)\setminus \mathfrak C(F)$. Then there exists a unique nonempty part~$H\in\P(\mathfrak C(F))$ such that $R\subset\O(H)$ and either~$R=\R(H)$, or~$R\cap \I(H)\ne\varnothing$.

\begin{proof}
If~$R\subset\O(F)$, then by lemmas~\ref{ll39} and~\ref{l41} the set~$R$ is a boundary of a nonempty part~$H\in\P(\mathfrak C(F))$. Hence $R\subset H\subset\O(H)$. 
Clearly, a boundary of a nonempty part of~$\P(\mathfrak C(F))$  neither  is a boundary  nor intersect the interior of another nonempty part of~$\P(\mathfrak C(F))$.

Let~$R\not\subset\O(F)$. Then  there exists a part~$G_{i,i+1}'\in\P(\mathfrak C(F))$ such that
$R\cap \I(G_{i,i+1}')\ne\varnothing$. Prove, that $R\subset\O(G_{i,i+1}')$. It is obvious in the case 
when~$R$ is independent with~$Q_{i,i+1}'$, since then~$R\subset G_{i,i+1}'$. Hence it is enough to consider
the case when these two  sets are dependent. In this case by lemma~\ref{ll39} the set~$R$ separates
one vertex~$x\in\O(F)$ from other vertices of~$\O(F)$.
Obviously, $x\in Q_{i,i+1}'=\{q_i,u_{i,i+1},q_{i+1}\}$. There are two possible cases.

\q1. If~$x=q_i$ (the case~$x=q_{i+1}$ is similar), then by lemma~\ref{ll37} the set~$R$ consists of two vertices
of the part~$G_{i,i+1}\subset\O(G_{i,i+1}')$ and a vertex~$q_{i-1}\in\O(G_{i,i+1}')$.

\q2. Otherwise, $x=u_{i,i+1}$. In this case $u_{i,i+1}\ne p$. Since~$p\in\O(F)$ and~$p$ is adjacent to~$u_{i,i+1}$, then~$p\in R$. Thus~$R$ is independent with~$Q_{i,i+1}$ (otherwise by lemma~\ref{l300} the flower~$F$ is not maximal). Hence, $R\subset G_{i,i+1}\subset\O(G_{i,i+1}')$.

It remains to notice, that the neighborhood of another part of~$\P(\mathfrak C(F))$ does not intersect~$\I(G_{i,i+1})$ and, consequently, does not contain~$R$.
\end{proof}
\end{thm}

\subsection{A complex of big cut}

\begin{defin}
1) We call a nontrivial cut~$M\in \mathfrak M_3$ {\it big}, if~$V(M)$ is not a subset of the neighborhood of any 
triple cut or nondegenerate flower. 

2) Define the {\it complex} $\mathfrak C(M)$ of a big cut~$M$ as a set of all 3-cutsets contained in~$V(M)$, 
except  boundaries of this cut~$T^M_1$ and~ $T^M_2$, which we call {\it boundaries} of~ $\mathfrak C(M)$.
\end{defin}

By corollary~\ref{l31c1} every set~$R\in \mathfrak C(M)$ is contained in the cut~$M$ 
(i.e., contains a vertex of each edge of the cut~$M$). All such sets, except boundaries of~$M$, belong 
to~$\mathfrak C(M)$.
As we know,~$\P(\mathfrak C(M))$ consists of normal parts~$G^M_1$ and~$G^M_2$, and small parts~$\{x_1,x_2\}$ 
where~$x_1x_2\in M$. We set, that the neighborhood of~$G^M_i$ as a part of~$\P(\mathfrak C(M))$ is 
its neighborhood as a part of~$\P(M)$.

It follows from lemma~\ref{ll34}, that for any set $R\in\mathfrak R_3(G)\setminus \mathfrak C(M)$ 
there exists a unique nonempty part~$A\in \P(\mathfrak C(M))$ such that~$R\subset \O(A)$ and  either~$R=\R(A)$,
or~$R\cap\I(A)\ne\varnothing$.

Let $M=\{a_1a_2,b_1b_2,c_1c_2\}$ be a big cut. As it was proved above, there exist six four-petal
flowers on vertices of~$V(M)$ (they are~$(b_1;a_1,a_2,c_2,c_1)$ and five similar flowers). 
Since the cut~$M$ is not contained in the neighborhood of a nondegenerate flower, all these flowers are
maximal and, of course, degenerate.

\subsection{Small complexes}

\begin{defin}
Let  triple complexes, complexes of nondegenerate flower and complexes of big cut
be  {\it big} complexes. All cutsets not belonging to any big complex we shall divide into complexes of
one or two cutsets. We call such complexes {\it small}.

Let the vertex set of each  (big or small) complex $\cal C$ be the union $V(\cal C)$ of all cutsets of~$\cal C$.
 \end{defin}

Each single  cutset form a single  complex (which is small). Further we describe  other small complexes.

Let $T\in\mathfrak R_3(G)$ be a nonsingle cutset,  not belonging to any big complex. Note, that
then~$|\P(T)|=2$. Indeed, otherwise by lemma~\ref{l3v0} any  3-cutset dependent with~$T$ is subordinated 
to~$T$, i.e.~$T$ is a line of triple complex. In addition, if cutset~$S\in\mathfrak R_3(G)$ is dependent 
with~$T$, then~$|\P(S)|=2$ (otherwise~$T$ is subordinated to~$S$ and belongs to a triple complex with line~$S$)
and $T\cap S=\varnothing$ (otherwise the sets~$T$ and~$S$ generate a flower).

\begin{lem}
\label{lmk0}
Let a cutset~$T=\{x,y,z\}$ be such that~$|\P(T)|=2$ and~$T$ can be complemented by each of edges~$xx_1$ and~$yy_1$, 
lying in different parts of~$\P(T)$. Then~$T$ can be complemented by both these edges simultaneously
(i.\,e. $\{xx_1,y_1y,z\}\in\mathfrak M_2(G)$) if and only if the vertices~$x$ and~$y$ are not adjacent.

\begin{proof}
Let the vertex~$x_1$ lie in a connected component~$H$ of the graph~$G-T$. By lemma~\ref{ll31} we know, that~$x_1$ is the only vertex of the component~$H$ adjacent to~$x$. Further consider a cut~$M_y=\{x,y_1y,z\}\in\mathfrak M_1(G)$. Obviously, $H\cup\{y\}$ is a  connected component of the 
graph~${G-M_y}$. By lemma~\ref{ll31} the cut~$M_y$ can be complemented by an edge~$xx_1$ if and only if~$x_1$
is the only vertex of the set~$H\cup\{y\}$, adjacent to~$x$. The last fact is equivalent to that $x$ and $y$ are not adjacent.
\end{proof}
\end{lem}

\begin{lem}
\label{lmk1}
Let $T=\{x,y,z\}$ be a nonsingle cutset not belonging to any big complex. Then the following statements hold. 

$1)$ For any cutset~$S\in\mathfrak R_3(G)$ dependent with~$T$ exactly one part of~$\P(\{S,T\})$ is small.  Vertices of this part form a singular edge, both cutsets~$S$ and $T$ can be complemented by this edge.

$2)$ All edges which complement the cutset~$T$ lie in the same part of~$\P(T)$.  
Moreover, the set~$T$ can be complemented by all these edges simultaneously.

\begin{proof}
1) As it was proved above, $T\cap S=\varnothing$, i.e. by corollary~\ref{l1c2} at least one  part of~$\P(\{S,T\})$ is small. From~$\P(T)=\P(S)=2$ it follows, that there is exactly one small part. By theorem~\ref{t30} 
we know, that  vertices of this part form a singular edge, by which both cutsets $S$ and $T$ can be complemented.

2) At first notice, that the cutset~$T$ cannot be complemented by two edges, lying in the different parts of~$\P(T)$
simultaneously. Indeed, otherwise~$T$ is an inner set of a cut of~$\mathfrak M_2(G)$, i.e. belongs to a big complex.

Now suppose, that the cutset~$T$ can be complemented by each of edges~$xx_1$ and $yy_1$, and these two edges lie in different parts of~$\P(T)$. Then by lemma~\ref{lmk0} the vertices~$x$ and~$y$ are adjacent. Consider a set~$S\in\mathfrak R_3(G)$ dependent with~$T$ (such a set exists, since~$T$ is nonsingle). We know, 
that~$T\cap S=\varnothing$, hence, $S$ separates~$z$ from~$\{x,y\}$. Then from item 1 of this lemma it follows, that the set~$T$ can be complemented by a singular edge~$zz_1$. Without loss of generality assume, that~$y_1$ and~$z_1$ are in different connected components of the graph~$G-T$. Then, similarly to proved above, the vertices~$y$ 
and~$z$ are adjacent and the set~$S$ cannot split~$T$, we obtain a contradiction.

If the set~$T$ can be complemented by edges~$xx_1$ and $xx_2$, then by lemma~\ref{ll31} each of the vertices
$x_1$ and $x_2$ is the only vertex adjacent to~$x$ in the connected component of the graph~$G-T$ containing this vertex. Since $d(x)\ge3$, it follows, that  $x$ is adjacent to a vertex~$y\in T$ and similarly to written above we obtain a contradiction.

Thus, all edges, by which we can complement the set~$T$ are in the same part of~$\P(T)$. Then by lemma~\ref{ll31} it follows, that the cutset~$T$ can be complemented by all these edges simultaneously.
\end{proof}
\end{lem}

Now we can describe all types of small complexes, parts of decomposition of the graph by a small complex and neighborhoods of these parts.

\begin{defin}
For any cut~$M=\{x_1x_2,y,z\}$ which boundaries  are nonsingle cutsets not belonging to any complex
we define the complex~$\mathfrak C(M)$ as a set consisting of both boundaries of~$M$. We call the cut~$M$ 
{\it small}, and complex~$\mathfrak C(M)$ --- a {\it complex of small cut}.

All other cutsets, not belonging to  big complexes or complexes of small cuts, form complexes consisting 
of one cutset.
\end{defin}

It is easy to see, that the complex of small cut~$M=\{x_1x_2,y,z\}$ splits $G$ into three parts:
$G^M_1$, $G^M_2$ and~$\{x_1,x_2,y,z\}$. All these parts are normal. But the part~$\{x_1,x_2,y,z\}$ 
is splitted by cutsets dependent with  boundaries of $M$,  and there is a  small 
part~$\{x_1,x_2\}\in \P({\mathfrak R_3}(G))$. This part is also empty. We set, that neighborhoods of~$G^M_i$ as a  part of~$\P(\mathfrak C(M))$ and as a part of $\P(M)$ coincide.

For every small complex ${\cal C}=\{T\}$ let us define the neighborhood of a part of $\P({\cal C})$. If~$T$ 
is a single set, then for every part~$H\in\P({\cal C})$ we set~$\O(H)=H$. Otherwise, let~$\P({\cal C})=\{H_1,H_2\}$.  By lemma~\ref{lmk1}, the ends of  all edges  which complement the set~$T$ are in the same part of~$\P({\cal C})$ --- let it be~$H_1$. 
We set~$\O(H_1)=H_1$. 
Let us define neighborhood of the other part~$H_2$. We complement the cutset $T$ to a maximal cut~$M$.
Clearly,~$H_2\in\P(M)$. Let the neighborhood of~$H_2$ as a part of~$\P({\cal C})$ be its neighborhood as a part of~$\P(M)$.

It follows from lemma~\ref{ll34}, that for any small complex~${\cal C}$ and any cutset~$R\in\mathfrak R_3(G)\setminus {\cal C}$ there exists a unique nonempty part~$A\in \P({\cal C})$ such that $R\subset \O(A)$ and $R\cap\I(A)\ne\varnothing$.

Let us describe all small complexes in details.

\begin{lem}
\label{lmk2}
Let $T=\{x,y,z\}$ be a nonsingle cutset not belonging to any big complex. Then at least one of the following three statements holds.

$1^\circ$ The cutset~$T$ is trivial.

$2^\circ$ The cutset~$T$ is a boundary of a big complex. All cutsets of this complex and all edges  which complement the set~$T$ lie in the same part of~$\P(T)$.

$3^\circ$ Exactly one edge~$xx_1$ complements the cutset~$T$, this edge is singular, and each cutset
dependent with~$T$ contains~$x_1$  and separates~$x$ from  $\{y,z\}$.

\begin{proof}
Since the cutset~$T$ is nonsingle, then there exists a cutset $S\in\mathfrak R_3(G)$ dependent with~$T$.
By lemma~\ref{lmk1} there is one small part in~$\P(\{S,T\})$ and its two vertices are ends of a singular edge, which complements the cutsets~$S$ and~$T$. Without loss of generality we may assume, that this edge is~$xx_1$.

Further we consider several cases.

\q1. {\it Let~$xx_1$ be the only edge which  complements~$T$.} Then for every set~$R\in\mathfrak R_3(G)$ dependent with~$T$ we have~$\{x,x_1\}\in\P(\{R,T\})$. Hence, $x_1\in R$ and $R$ separates~$x$ from $\{y,z\}$. Thus in this case statement~$3^\circ$ holds.

\q2. {\it Let  an edge~$yx_1$ also complement the cutset~$T$.} Denote by~$H$ a connected component of the graph~$G-T$ containing~$x_1$.  By lemma~\ref{lmk1}, the set  $T$ can be complemented by edges~$yx_1$ and~$xx_1$ simultaneously. 
Then by item~2 of remark~\ref{nr1} we obtain, that~$H=\{x_1\}$ and~$T$ is a trivial cutset.
In this case statement~$1^\circ$ holds.

\q3. {\it Let  an edge~$yy_1$ also complement the cutset~$T$ (where $x_1\ne y_1$).} Then by lemma~\ref{lmk1} the cutset~$T$ can be complemented by the edges~$xx_1$ and~$yy_1$ simultaneously, i.e.
$M=\{xx_1,yy_1,z\}\in\mathfrak M_2(G)$. Inner sets of the cut~$M$ generate a flower~$(z;x,x_1,y_1,y)$, 
which is contained in some big complex~$\cal C$. 
Then~$T\subset V(\cal C)$. But by condition of lemma~$T \not\in\cal C$, hence, the cutset~$T$ is a boundary
of the complex~$\cal C$. Since a boundary of a complex cannot split its vertex set, then~$V(\cal C)$ is contained in the part of~$\P(T)$ which contains the vertices~$x_1$ and~$y_1$. 
In this case statement~$2^\circ$ holds.
\end{proof}
\end{lem}

\begin{rem}
\label{rmk1}
1) The statement~$1^\circ$ cannot be fulfilled simultaneously with one of statements~$2^\circ$ or~$3^\circ$.

2) It is easy to see from the prove of lemma~\ref{lmk2}, that if a nonsingle 3-cutset not belonging to any big complex can be complemented by exactly one edge, then statement~$3^\circ$ of lemma~\ref{lmk2} holds for this cutset.

3) A boundary of a triple complex or of a complex of big cut always can be complemented by an edge, lying in the part containing all vertices of this complex. A boundary~$Q_{i,i+1}'$ of a complex of  flower~$\mathfrak C(F)$ 
cannot be complemented by such edge in the only case: the flower~$F$ is nonsingular,
$u_{i,i+1}=p$ and both parts~$G_{i-1,i}$ and~$G_{i+1,i+2}$ are nonempty. It is easy to see, that in this case the set $Q_{i,i+1}'$ is single.

It follows from written above, that if a boundary of a big complex is not a single cutset and do not belong to another big complex, then it cannot be complemented by an edge lying in the part not containing all vertices of this complex.

In particular, all cutsets, belonging to a complex of small cut cannot be boundaries of big complexes. 
Clearly, such cutsets also cannot be trivial or single. Hence, each set of the complex of a small cut~$M=\{x_1x_2,y,z\}$ can be complemented by an edge~$x_1x_2$.
\end{rem}

\begin{lem}
\label{lmk3}
Let~$T=\{x,y,z\}$ be a nonsingle cutset, which do not belong to any big complex and is not a boundary of big
complex. Let~$T$ can be complemented by exactly one edge~$xx_1$. 
Then~$T_1=\{x_1,y,z\}$ is a cutset and one of  two following statements holds.

$1^\circ$ The cutset~$T_1$ is single.

$2^\circ$ Two cutsets  $T$ and~$T_1$ form a complex of small cut.

\begin{proof}
The set~$T_1$ can be not a cutset in the only case: if $\{x_1\}$ is a connected component of the graph~$G-T$. 
But then the cutset~$T$ is trivial and can be complemented by edges~$yx_1$ and~$zx_1$ too. We obtain a contradiction.

Thus,~$T_1$ is a cutset. Suppose, that it is nonsingle. 
We need to prove, that in this case the cutsets~$T$ and~$T_1$ form a complex of small cut, i.e.~$T_1$ 
does not belong to any big complex. 

Note, that by lemma~\ref{lmk2} the edge~$xx_1$ is singular and any  3-cutset~$S$ dependent with~$T$ 
contains~$x_1$ and separates~$x$ from $\{y,z\}$. Hence, $S$ cannot be dependent with~$T_1$. 
Thus, every  3-cutset dependent with~$T_1$ is independent with~$T$.

Now assume, that~$T_1$ belongs to a big complex~${\cal C}$. Then~$T_1$ must be dependent with at least one 
cutset~$S\in{\cal C}$. Since~$S$ is dependent with~$T_1$ and independent with~$T$, it must contain the vertex~$x$. But then  $T\subset V({\cal C})$, i.e. either~$T$ belongs to~${\cal C}$, or~$T$ is a boundary
of~${\cal C}$. We obtain a contradiction.
\end{proof}
\end{lem}

\begin{cor}
\label{lmk3c1}
Let a complex consist of one cutset. Then this cutset can be:

$1)$ a single cutset;

$2)$ a trivial cutset;

$3)$ a boundary of big complex;

$4)$ a cutset which can be complemented by exactly one edge, and the other boundary of resulting cut
is a single cutset.
\end{cor}

\section{Relative position of complexes}

In previous section we have described all types of complexes and for each type we have investigated some properties. Let us repeat  properties, that hold for all types of complexes.

For every complex~${\cal C}$ a boundary of any nonempty part~$A\in\P({\cal C})$ is a 
3-cutset, which do not split~$V({\cal C})$ (but can split~$A$). 
Let~$R=\R(A)$. If~${\cal C}$ consists of more than one cutset, then $\P(R)$ contains exactly one part which do not intersect~$\I(A)$. Denote this part by~$\overline{A}$. 
In this case the neighborhood of~$A$ is constructed as follows:
the cutset~$R$ is complemented   by all possible edges lying in~$\overline{A}$. 
Let $M$ be the resulting cut.
After that we set, that~$\O(A)$ is a neighborhood of~$A$ as a part of $\P(M)$. 
Note also, that if~$A\in\P({\cal C}_1)$ 
and~$A\in\P({\cal C}_2)$ (it is possible, for example, if~${\cal C}_1$ is a big complex and~${\cal C}_2=\{R\}$, where~$R$ is a bound of~${\cal C}_1$ and~$|\P(R)|=2$), then neighborhoods of~$A$ in both cases coincide.

\begin{defin}
For any complex ${\cal C}$ and any cutset~$T\in\mathfrak R_3(G)\setminus {\cal C}$ 
we say, that $T$ {\it belongs} to a nonempty part $A\in \P({\cal C})$, 
if $T\subset \O(A)$ and either~$T=\R(A)$, or $T\cap\I(A)\ne\varnothing$.
\end{defin}

In previous section it was proved, that for any complex~${\cal C}$  any 
cutset~$T\in\mathfrak R_3(G)\setminus {\cal C}$ belongs to exactly one nonempty part of~$\P({\cal C})$.

Our first aim is to show, that two cutsets belonging to one complex cannot belong to different parts of decomposition of the graph by another complex.
For this purpose we need the following lemmas.

\begin{lem}
\label{lts00}
Let ${\cal C}$ be a complex and $T\in\mathfrak R_3(G)\setminus {\cal C}$ be a set splitting~$V({\cal C})$.
Let~$T$ belong to a part~$A\in\P({\cal C})$ and~$R=\R(A)$. Then the following statements hold.

$1)$ The cutsets~$R$ and~$T$ are dependent.

$2)$ The cutset~$T$ separates exactly one vertex~$x\in R$ from other vertices of~$V({\cal C})$.

$3)$ The cutset~$T$ consists of two vertices of the part~$A$ and a vertex~$y\not\in A$ such that both cutsets~$R$ and~$T$ can be complemented by the edge~$xy$.

\begin{proof}
1) Note, that $T\not\subset A$, since otherwise  $T$ cannot split~$V({\cal C})$. 
Hence, $T\ne R$. Since~$T$ belongs to the part~$A$ we obtain, that $T\cap\I(A)\ne\varnothing$.
Thus, $R$ splits~$T$, consequently, these sets are dependent.

2) It follows from previous lemmas, that the cutset~$T$ separates exactly one vertex~$x\in V({\cal C})$ from other vertices of the set~$V({\cal C})$: for  a complex of big or small cut it follows from lemma~\ref{ll34}, 
for a triple complex --- from lemma~\ref{ll35}, for a complex of flower --- from lemma~\ref{ll39}. 
In the case~$|{\cal C}|=1$  this statement is obvious. Since~$T$ and~$R$ are dependent, then~$x\in R$.

3) Add the cutset~$R$ to a cut~$M$ by all possible edges lying in the part~$\overline{A}$
(remind, that~$\overline{A}$ is the only part of $\P(R)$ not intersecting~$\I(A)$). Let~$M'$ be a maximal
cut containing~$M$. Since~$T\subset\O(A)=A\cup V(M)$, the cutset~$T$ cannot generate a flower with both boundaries of the cut~$M'$. If~$T\not\subset V(M')$, then statement 3 of our lemma follows from item~2 of lemma~\ref{ll34} for the cut~$M'$ and the cutset~$T$.  If~$T\subset V(M')$, then we have~$M'\ne M$. It is possible only if~${\cal C}$ is a triple complex or a complex of flower, $M\in\mathfrak M_1(G)$, and~$M'\in\mathfrak M_2(G)$. In this case the statement we prove is clear.
\end{proof}
\end{lem}

\begin{cor}
\label{lts00c1}
Let complexes~${\cal C}_1$ and~${\cal C}_2$ be such that ${\cal C}_2=\{T\}$ and the cutset~$T$ 
splits~$V({\cal C}_1)$. Let the cutset $T$ belong to a part~$A\in\P({\cal C}_1)$, $R=\R(A)$ and~$\overline{A}$ is a part of~$\P(R)$ not intersecting~$\I(A)$.  Then~$|\P(R)|=|\P(T)|=2$. In addition, parts of~$\P(T)$ can be denoted by~$B$ and~$\overline{B}$ such that the  following statements hold.

$1)$ $|\overline{A}\cap\overline{B}|=2$.

$2)$ $\O(\overline{B})=\overline{B}\subset\O(A)$.

$3)$ All cutsets of the complex~${\cal C}_1$ belong to the part~$B\in\P({\cal C}_2)$.

\begin{proof}
By lemma~\ref{lts00} the cutsets~$R$ and~$T$ are dependent, $T$ separates a vertex $x\in R$
from  other vertices of the set~$V({\cal C}_1)$ and both cutsets~$T$ and~$R$ can be complemented by an edge~$xy$ 
(where $y\in T$). Moreover, $T\setminus A=\{y\}$.
Since~$T$ does not belong to big complexes, $T\cap R=\varnothing$ and $|\P(R)|=|\P(T)|=2$.
Let~$\P(T)=\{B,\overline{B}\}$ where~$x\in\overline{B}$. Let us check that all the statements hold.

1) Since~$T\cap R=\varnothing$ we have $T\cap\overline{A}=\{y\}$ and~$R\cap\overline{B}=\{x\}$.
By  corollary~\ref{l1c2} that  means~$\overline{A}\cap\overline{B}=\{x,y\}$.

2) By lemma~\ref{lmk1} all edges  which  complement~$T$ lie in the part~$\overline{B}$, 
hence, $\O(\overline{B})=\overline{B}$. In addition, 
$\overline{B}\setminus A=\overline{B}\cap\I(\overline{A})=\{y\}\subset\O(A)$, thus,
$\overline{B}\subset\O(A)$.

3) Since $\O(\overline{B})=\overline{B}$, then only cutsets contained in~$\overline{B}$ can belong to the part~$\overline{B}$. However, $V({\cal C}_1)\cap\overline{B}\subset\{x,y\}$,
hence,  cutsets of the complex~${\cal C}_1$ cannot belong to the part~$\overline{B}$.
\end{proof}
\end{cor}

\begin{lem}
\label{lts0}
For any maximal nontrivial cut~$M$ one of two following  statements holds.

$1^\circ$ All $3$-cutsets contained in~$M$ (i.e. inner sets and boundaries of~$M$)
belongs to some complex~${\cal C}$.

$2^\circ$ A vertex set of any complex~${\cal C}$ is contained in a neighborhood of some part of~$\P(M)$.

\begin{proof}
By lemma~\ref{ll34} if a cutset~$T\in\mathfrak R_3(G)$ is not contained in the neighborhood of any part of~$\P(M)$, then $T$ with both boundaries of~$M$ generates a flower,  which is  contained in a maximal flower~$F$. Then all inner sets and boundaries of the cut~$M$ belong to the complex of flower~$F$, and statement~$1^\circ$ holds.

Let any cutset~$T\in \mathfrak R_3(G)$ be contained in a neighborhood of some part of~$\P(M)$. We shall prove, that then statement~$2^\circ$ holds.
In the case~$|{\cal C}|=1$ it is clear. Assume, that~$|{\cal C}|>1$ and consider several cases.

\q{a}.  {\it Let~$\cal C$ be a complex of big or small cut.} Note, that the vertex set of another cut~$M'$ 
is contained in the neighborhood of some part of~$\P(M)$.  It is clear, since any two vertices of~$V(M')$
are either adjacent, or belong to a 3-cutset contained in~$M'$~--- in both cases they cannot lie in interiors of different parts of~$\P(M)$. Hence  statement~$2^\circ$  for a complex of big or small cut immediately follows.

\q{b}. {\it  Let~${\cal C}$ be a triple complex with line~$S$.} Since the cut~$M$ is nontrivial, then~$S$ is independent with both boundaries of~$M$. Let~$S\subset G_1^M$. Then a vertex set of any cut with boundary~$S$ is contained in $\O(G_1^M)$, consequently, $V({\cal C})\subset\O(G_1^M)$.

\q{c}. {\it It remains to consider the case when~${\cal C}$ is the  complex of a flower $F$.} In this case~$V(F)$ 
is contained in the neighborhood of some part of~$\P(M)$ (let it be~$G_1^M$), since any two vertices of the 
set~$V(F)$ are either adjacent, or belong to a 3-cutset.
Hence, any set of the flower~$F$ is independent with~$T_2^M$, consequently, $G_2^M$ is contained 
in some part of~$\P(F)$  --- let it be~$G_{i,i+1}$. Thus, only~$u_{i,i+1}$ can be a vertex of~$V(\cal C)$ 
not lying in~$\O(G_1^M)$ (and, hence, contained in~$\I(G_2^M)$). It is possible when~$u_{i,i+1}\ne p$. 
In this case~$u_{i,i+1}$  is the only vertex of the set~$G_{i,i+1}$ (and, consequently, of the 
set~$G_2^M\subset G_{i,i+1}$) which is  adjacent to~$p$. 
Since $u_{i,i+1}\in\I(G_2^M)$, we have, that $p\in T_2^M$. Thus $p$ cannot be an end of edge~$px\in M$: otherwise~$p$ is adjacent to~$u_{i,i+1}$,  $x$ and only them, that is impossible.
By lemma~\ref{ll31}, hence, the cut~$M$ can be complemented by an edge~$pu_{i,i+1}$. We obtain a contradiction with
maximality of~$M$.  Consequently,~$V({\cal C})\subset \O(G_1^M)$.
\end{proof}
\end{lem}

\begin{lem}
\label{lts1}
Let ${\cal C}_1$ and~${\cal C}_2$ be two complexes. Then all cutsets of~${{\cal C}_2\setminus{\cal C}_1}$ 
belong to one part~$A\in \P({\cal C}_1)$.  Moreover,~$V({\cal C}_2)\subset\O(A)$.

\begin{proof}
In the case~$|{\cal C}_2|=1$ the statement of lemma is obvious. In the case~$|{\cal C}_1|=1$ this statement
immediately follows from corollary~\ref{lts00c1}. Thus it is enough to consider the case~$|{\cal C}_1|>1$ 
and~$|{\cal C}_2|>1$. Hence we obtain that the sets~$V({\cal C}_1)$ and~$V({\cal C}_2)$ does not contain each other. (The vertex set of a big complex by construction cannot be a subset of the vertex set of another complex.
For a small complex it is possible only if it consists of one cutset, which is a boundary of big complex.)
Consider a vertex~$u\in V({\cal C}_2) \setminus V({\cal C}_1)$ and a part~$A\in\P({\cal C}_1)$ containing~$u$.
Since $u\not\in V({\cal C}_1)$, we obtain, that~$u\in\I(A)$.

Note, that since~$|{\cal C}_1|>1$, the neighborhood of each part of~$B\in \P({\cal C}_1)$ is contained in~$V({\cal C}_1)\cup B$. Hence no  3-cutset can intersect interiors of two parts of~$\P({\cal C}_1)$, consequently, all cutsets of the complex~${\cal C}_2$ intersecting~$\I(A)$ belong to~$A$. 

Let us prove, that $V({\cal C}_2)\subset\O(A)$. Consider several cases.

\q1. {\it Let~${\cal C}_2$ be a complex of small cut.} If both cutsets of this complex contain~$u$, then these cutsets belong to~$A$ and, hence, are contained in~$\O(A)$.
If only one cutset~$T\in{\cal C}_2$ contains~$u$, then $V({\cal C}_2)$ consists of vertices of the set~$T$ and a vertex~$u_1$ adjacent to~$u$. Since~$u\in\I(A)$ we obtain~$u_1\in A$.

\q2. 
{\it Let ${\cal C}_2$ be a big complex. }

\q{2.1}. {\it Let~${\cal C}_1$ be a complex of big or small cut.}
 Then the statement we prove follows from
lemma~\ref{lts0}. 

\q{2.2}. {\it Let ${\cal C}_1$ be the triple complex with line~$S$.}
 In this case the cutset~$S$ is independent with all cutsets
of~${\cal C}_2$, hence, $V({\cal C}_2)$ is contained in some part of~$\P(S)$,
and every part of~$\P(S)$ is a neighborhood of correspondent part of~$\P({\cal C}_1)$.

\q{2.3}. {\it Let ${\cal C}_1$ be the complex of a flower~$F$.}
 Then $A=G'_{i,i+1}$ for some~$i$.
Consider the set~$M'_{i,i+1}$. If~$M'_{i,i+1}$ is a maximal cut, the statement we prove immediately  follows from lemma~\ref{lts0}. This statement is also clear if~$M'_{i,i+1}$ contains no edge (in this case it follows from lemma~\ref{ll36} and lemma~\ref{ll37}, that~$M'_{i,i+1}$ is a single cutset, which is the boundary of the part~$A$). Hence it is enough to consider the case, when $M'_{i,i+1}$ contains at least one edge, but is not a maximal cut. 
It is clear from definition~\ref{dor}, that it is possible only when $F$ is a nonsingular flower, exactly one of the parts~$G_{i-1,i}$ and~$G_{i+1,i+2}$ is empty (without loss of generality assume, that it is~$G_{i+1,i+2}$) and the center~$p$ is adjacent to exactly one vertex of the part~$G_{i,i+1}$ (denote this vertex by~$u'$). Consider this case in details.

It follows from lemma~\ref{lor0}, that the cut~$M'_{i,i+1}$ can be complemented only by the edge~$pu'$. 
Denote the resulting maximal cut by~$M$. By lemma~\ref{lts0} the set~$V({\cal C}_2)$ is contained in the neighborhood of some part of~$\P(M)$.  Note, that~$G_{i,i+2}=\O(A)$ is a neighborhood of one part of~$\P(M)$.
If $V({\cal C}_2) \subset G_{i,i+2}$, then the statement we prove is fulfilled.
The other part of~$\P(M)$ is~$G_{i+2,i}$, and $\O(G_{i+2,i})=G_{i+1,i}\cup\{u'\}$. 
Thus we can consider the case when~$V({\cal C}_2)\subset G_{i+1,i}\cup\{u'\}$ 
and~$V({\cal C}_2)\not\subset V(M)$. Note, that then~$V({\cal C}_2)\cap\I(A)=\{u'\}$ and~$u'=u$. 
Consider a cutset~$T\in{\cal C}_2$ containing~$u$. Since~${\cal C}_2$ is a big complex, there is a 
cutset~$R\in{\cal C}_2$ dependent with~$T$. Whence follows that~$T\ne\{q_i,u,q_{i+1}\}$. But then~$T$ is dependent with~$Q_{i,i+1}$, i.e.~$T$ splits~$V(F)$. 
By lemmas~\ref{ll36} and~\ref{ll37} we have, that~$T=\{q_i,u,q_{i+2}\}$, and, moreover,
$T$ is the only cutset of the complex~${\cal C}_2$ which contains~$u$. Since $R$ and~$T$ are
dependent, then $q_{i+1}\in R$, i.e. $\{u,q_i,q_{i+1},q_{i+2}\}\subset V({\cal C}_2)$. 
In addition, it follows from lemma~\ref{l31}, that~$q_{i+2}\not\in R$.

Since the vertex~$u$ belongs to exactly one cutset of~${\cal C}_2$, this complex is not a complex of big cut. 
It remains to consider two cases: ${\cal C}_2$ is a complex of flower or a triple complex.

\q{2.3.1}. {\it Let ${\cal C}_2$ be the complex of a flower~$F'$.}
 Clearly, in this case~$u$ is a petal of~$F'$ and the flower~$F'$ has four  petals. In addition, $V({\cal C}_2)=V(F')$, since otherwise ${\cal C}_2$ would be a complex 
of big cut. Note, that only~$q_i$ can be a center of $F'$ (since  another vertex cannot belong to both cutsets
$T, R \in {\cal C}_2 = \mathfrak R(F')$). But then a cutset~$R'=\{q_i,p,q_{i+1}\}$ contains the center of~$F'$ 
and is dependent with the cutset~$T \in \mathfrak R(F)$. Hence by lemma~\ref{l300}  we have, that~$R'\in \mathfrak R(F)$ and, since~$F'$ has only  4 petals, then $R'=R$ and $F'=(q_i;u,q_{i+1},q_{i+2},p)$. 
Thus $V({\cal C}_2)\subset\O(A)$.

\q{2.3.2}. {\it Let~${\cal C}_2={\cal C}(N)$ be a triple complex, and a triple cut~$N=M_1\cup M_2\cup M_3$ has line~$S$.}
The line of a triple complex splits the graph into three parts, hence by item~2  of lemma~\ref{ll39} it cannot split  vertex set of the flower~$F$. Then~$S\ne T$. Moreover, $S$ and $T$ are independent, since otherwise
three vertices of the set~$T$ would be in three different connected components of the graph~$G-S$, i.e.
the cutset~$S$ splits~$\{q_i,q_{i+2}\}\subset V(F)$. Now without loss of generality we may assume, that the cutset~$T$ is contained in the cut~$M'_1$. Remind, that the vertex~$u\in T$ do not belong to other cutsets
of the complex~${\cal C}_2$. In particular, $u\notin S$. Thus there exists an edge~$ux\in M_1'$, 
but in this case the vertex~$u$ belongs to at least two cutsets of the complex~${\cal C}_2$. We obtain a  contradiction and show, that this case  is impossible.

\smallskip

Thus we have proved, that in all cases~$V({\cal C}_2)\subset\O(A)$. Hence any cutset of the complex~${\cal C}_2$ either intersects~$\I(A)$ and, consequently, belongs to~$A$, or is contained in~$V({\cal C}_1)$. 
Note, that a  3-cutset~$R\subset V({\cal C}_1)$ can either belong to the complex~${\cal C}_1$, or be a boundary
of some part of~$\P({\cal C}_1)$. Let a cutset~$R\in{\cal C}_2$ be a boundary of a part~$A_1\in\P({\cal C}_1)$ different from~$A$. Then~$R$ is independent with all sets of~${\cal C}_2$, that is possible only 
when~${\cal C}_2$ is a complex of small cut. In this case~$V({\cal C}_2)\subset \O(A_1)$, hence $V({\cal C}_2)\cap\I(A)=\varnothing$. We obtain a contradiction.
\end{proof}

\end{lem}

\begin{rem}
Note, that two complexes ${\cal C}_1$ and ${\cal C}_2$ can have nonempty intersection. For example,  complexes of two nondegenerate flowers can have common boundary cut.
\end{rem}

\goodbreak
\begin{defin}
Let $\mathfrak C=\{{\cal C}_1,\ldots,{\cal C}_n\}$ be the set of all complexes of the graph~$G$. 
Denote by~$A_{i\supset j}$ the part of~$\P({\cal C}_i)$ to which all cutset of the complex~${\cal C}_j$  belong. 
Let us say, that the complex~${\cal C}_j$ belongs to the part~$A_{i\supset j}$.

For each complex~${\cal C}_i$ we denote by~$\mathfrak C_i$ the decomposition of all other complexes into classes: two complexes~${\cal C}_j$ and~${\cal C}_\ell$ belongs to the same class of this decomposition if and only if~$A_{i\supset j}=A_{i\supset \ell}$.

We say that a complex ${\cal C}_i$ {\it separates} ${\cal C}_j$ from~${\cal C}_\ell$,  if they belong  to different classes of the decomposition~$\mathfrak C_i$. We call complexes~${\cal C}_i$ and~${\cal C}_j$ {\it neighboring}, if no other complex separates ${\cal C}_i$ from~${\cal C}_j$.
Denote by~$T(G)$ a hypergraph which vertices are complexes of the graph~$G$, and hyperedges are all maximal
with respect to inclusion sets of pairwise neighboring complexes. We call the hypergraph~$T(G)$ a 
{\it hypergraph of decomposition} of~$G$.
\end{defin}

The construction of hypergraph of decomposition is described in details
in~\cite[section~2]{k05}. In this work we will use the following theorem (see~\cite[theorem~3]{k05}).

\begin{thm}
\label{tt}
Let for each element of the set~$V$   the decomposition of all other elements of~$V$
into classes  correspond.  Let for every $a,b,c\in V$ the following condition hold: if~$a$ separates $b$ from~$c$, 
then $b$ does not separate~$a$ from~$c$. Then the following statements hold.

$1)$ The hypergraph of this decomposition~$T(V)$ is a hypertree (i.e., any cycle of this hypergraph
is a subset of some hyperedge).

$2)$ Let for some vertex~$a\in V$ the hypergraph~$T(V)-a$ have connected components~$W_1,\ldots, W_\ell$. Then the element~$a$ decompose elements of~$V\setminus \{a\}$ exactly into classes~$W_1,\dots,W_\ell$.
\end{thm}

\begin{lem}
\label{ls1}
Let $B$  be a nonempty part of~$\P({\cal C}_j)$, different from~$A_{j\supset i}$. Then
$\O(B)\subset \O(A_{i\supset j})$. Moreover, if $B\not\subset A_{i\supset j}$, 
then  $|{\cal C}_j|=1$ and the only cutset of the complex~${\cal C}_j$ splits~$V({\cal C}_i)$.

\begin{proof}
Let~$R=\R(A_{i\supset j})$, $S=\R(A_{j\supset i})$, $T=\R(B)$. 

Note, that if $B\subset A_{i\supset j}$, then $\O(B)\subset \O(A_{i\supset j})$. 
Indeed, if  $y\in\O(B)\setminus B$, then the cutset~$T$ can be complemented by an edge~$xy$ (where $x\in T$).
If in addition~$y\not\in A_{i\supset j}$, then $x\in R$, and it follows from lemma~\ref{ll31}, that the cutset~$R$ can be also complemented by the edge~$xy$. Hence $y\in\O(A_{i\supset j})$.

Thus, it is enough to prove, that~$B\subset A_{i\supset j}$. We shall do it in all cases, except~1.3.
Consider several cases.

\q1. {\it Let ${\cal C}_j=\{T\}$.}  Clearly, then $S=T$.
This case is divided into the following subcases.

\q{1.1}. $T=R$. Then $T$ does not split~$V({\cal C}_i)$, hence, $V({\cal C}_i)\subset A_{j\supset i}$. 
Moreover, in this case ${\cal C}_i$ is a big complex, consequently, $A_{i\supset j}$  is a union
of all parts of~$\P(T)$ except~$A_{j\supset i}$. Thus, $B\subset A_{i\supset j}$.

\q{1.2}. {\it $T\ne R$ and~$T$ does not split $V({\cal C}_i)$. }
Then~$R\subset V({\cal C}_i)\subset A_{j\supset i}$. Hence, $R\cap\I(B)=\varnothing$, i.e.
by lemma~\ref{lds1} we have, that $R$ does not split $B$. Then~$B$ is contained in some part~$H\in\P(R)$. In addition~$H\subset A_{i\supset j}$, since~$T\subset A_{i\supset j}$ and $T\ne R$.
Thus, $B\subset A_{i\supset j}$.

\q{1.3}. {\it $T$ splits $V({\cal C}_i)$.}
In this case by corollary~\ref{lts00c1} we have, that $\P(T)=\{A_{j\supset i},B\}$ and 
$\O(B)=B\subset \O(A_{i\supset j})$.
Note, that it is the only case when~$B\not\subset A_{i\supset j}$.

\q2. {\it Let $|{\cal C}_j|>1$.}

Then $V({\cal C}_j)\not\subset V({\cal C}_i)$, consequently, 
$V({\cal C}_j)\cap \I(A_{i\supset j})\ne\varnothing$.
Moreover, in this case~$V({\cal C}_i)\subset\O(A_{j\supset i})\subset A_{j\supset i}\cup V({\cal C}_j)$.
Note, that the cutset~$T$ does not split~$A_{j\supset i}\cup V({\cal C}_j)$, consequently,
$A_{j\supset i}\cup V({\cal C}_j)$ is contained in some part~$H\in\P(T)$. 
In addition $H\ne B$, since~$A_{j\supset i}\subset H$. Then~$V({\cal C}_i)\cap \I(B)=\varnothing$ and
$B$ is contained in some part of~$\P({\cal C}_i)$.  We want to prove  that $B\subset A_{i\supset j}$.
For this purpose it is enough  to check, that  $T\cap\I(A_{i\supset j})\ne\varnothing$.

Suppose the converse. Then~$T$ either contains an inner vertex of another part of~$\P({\cal C}_i)$,
or~$T$ is a subset of~$V({\cal C}_i)$.
Consider these two cases in details. Let $\overline{A_{j\supset i}}$ be the part of~$\P(S)$ which contains~$V({\cal C}_j)$. 

\q{2.1}. {\it Let  a part~$F\in\P({\cal C}_i)$ such that $F\ne A_{i\supset j}$ and~$T\cap \I(F)\ne\varnothing$  exist. }
It is possible only if~${\cal C}_i=\{R\}$ and~$R$ splits~$V({\cal C}_j)$.
Then by corollary~\ref{lts00c1} we obtain, that~$|F\cap\overline{A_{j\supset i}}|=2$.
On the other side, $T\subset V({\cal C}_j)\subset\overline{A_{j\supset i}}$, hence,
$T\not\subset F$. But then~$T$ and~$R$ are dependent, that is impossible, since~$R\subset H\in\P(T)$.
Thus, this case is impossible.

\q{2.2}. {\it Let $T\subset V({\cal C}_i)$.} Then~$T\subset\O(A_{j\supset i})$. As we know,  $S=\R(A_{j\supset i})$ and~$S$ is independent with~$T$, consequently, there exists a cut~$M$ with boundaries~$S$ and~$T$. 
In addition, $V({\cal C}_j)\subset\overline{A_{j\supset i}}\cap H=V(M)$. Since $|{\cal C}_j|>1$, hence by the definition of complex it follows, that $V({\cal C}_j)=V(M)$ and~$M$ is a maximal cut. 

Since $V({\cal C}_i)\subset H\in\P(T)$, then the cutset~$T$ does not split~$V({\cal C}_i)$. 
Then it follows from~$T\subset V({\cal C}_i)$, that~$T$ is a boundary of the complex~${\cal C}_i$. 
Also note, that~${\cal C}_i\ne\{T\}$, since the cutset $T$ belongs to different from~$A_{j\subset i}$ part $B\in\P({\cal C}_j)$. Hence,~$B\in\P({\cal C}_i)$ and~$\O(B)=B\cup V(M)$.
Consequently, $V(M)\subset V({\cal C}_i)$ and we obtain a contradiction 
with~$V({\cal C}_j)\cap \I(A_{i\supset j})\ne\varnothing$.
Thus, this case is also impossible.
\end{proof}
\end{lem}

\begin{thm}
\label{ts2}
$1)$ The hypergraph~$T(G)$ is a  {\it hypertree} (i.e., each cycle of~$T(G)$ is a subset of some hyperedge).

$2)$ Let  ${\cal C}_i\in\mathfrak C$ and~$H_1,\ldots, H_\ell$ be connected components of the hypergraph~$T(G)-{\cal C}_i$. Then~$\mathfrak C_i=\{H_1,\ldots, H_\ell\}$.

\begin{proof}
Both statements of this theorem immediately follow from theorem~\ref{tt}, hence it is enough to verify the conditions of this theorem.

Suppose, that the complex~${\cal C}_i$ separates ${\cal C}_j$ from~${\cal C}_\ell$, i.e.
$A_{i\supset j}\ne A_{i\supset \ell}$. We need to prove, that~${\cal C}_j$ does not separate~${\cal C}_i$ 
from~${\cal C}_\ell$, i.e.~$A_{j\supset i}=A_{j\supset \ell}$.
By lemma~\ref{ls1} we have $\O(A_{i\supset \ell})\subset\O(A_{j\supset i})$, consequently,
$V({\cal C}_\ell)\subset\O(A_{j\supset i})$. On the other side, 
$V({\cal C}_\ell)\subset\O(A_{j\supset\ell})$.

Let~$A_{j\supset i}\ne A_{j\supset \ell}$. 
Then~$V({\cal C}_\ell)\subset\O(A_{j\supset i})\cap\O(A_{j\supset\ell})$.
In addition, $\O(A_{j\supset i})\ne A_{j\supset i}$, since otherwise
$V({\cal C}_\ell)\subset A_{j\supset i}$ and, consequently, the complex~${\cal C}_\ell$ belongs to the 
part~$A_{j\supset i}$.
Further we consider the following two cases.

\q1. {\it Let $|{\cal C}_j|>1$.} Then 
$V({\cal C}_\ell)\subset\O(A_{j\supset i})\cap\O(A_{j\supset\ell})\subset V({\cal C}_j)$.
Hence, ${\cal C}_j$ is a big complex and~${\cal C}_\ell=\{T\}$  where~$T=\R(A_{j\supset\ell})$. 
That is, $V({\cal C}_\ell)\subset A_{j\supset\ell}$.  Moreover, by lemma~\ref{ls1} in this 
case~$A_{j\supset \ell}\subset A_{i\supset j}$, hence, $V({\cal C}_\ell)\subset A_{i\supset j}$. 
But then the complex~${\cal C}_\ell$ belongs to the part~$A_{i\supset j}$, i.e. $A_{i\supset j}=A_{i\supset \ell}$. We obtain a contradiction.

\q2. {\it Let ${\cal C}_j=\{R\}$.} Since $\O(A_{j\supset i})\ne A_{j\supset i}$, then
the cutset~$R$ is nonsingle. Hence, $\P({\cal C}_j)=\{A_{j\supset i},A_{j\supset\ell}\}$ and by 
lemma~\ref{lmk1} we obtain, that all edges  which complement the cutset~$R$ lie in the part~$A_{j\supset\ell}$. Thus $\O(A_{j\supset\ell})=A_{j\supset\ell}$, i.e.
$V({\cal C}_\ell)\subset A_{j\supset\ell}$.
Hence~$R$ splits~$V({\cal C}_i)$, since otherwise by lemma~\ref{ls1} we 
have~$A_{j\supset\ell}\subset A_{i\supset j}$ and, consequently, the complex ${\cal C}_\ell$ 
belongs to the part~$A_{i\supset j}$, that contradicts the assumption. Let $S=\R(A_{i\supset j})$.
Since~$R$ splits~$V({\cal C}_i)$, then by corollary~\ref{lts00c1} we have~$\P(S)=\{A_{i\supset j},\overline{A_{i\supset j}}\}$ and~$|A_{j\supset\ell}\cap\overline{A_{i\supset j}}|=2$. However, $|V({\cal C}_\ell)|\ge 3$, consequently,  $V({\cal C}_\ell)\cap\I(A_{i\supset j})\ne\varnothing$ and the 
complex~${\cal C}_\ell$ belongs to the part~$A_{i\supset j}$. We have a contradiction.
\end{proof}
\end{thm}

\smallskip
Translated by D.\,V.\,Karpov.

\end{document}